\documentclass[onefignum,onetabnum]{siamart220329}
\allowdisplaybreaks

\usepackage{lipsum}
\usepackage[export]{adjustbox}
\usepackage{amsfonts}
\usepackage{graphicx}
\usepackage{epstopdf}
\usepackage{algorithmic}
\usepackage{cancel}
\usepackage{caption}
\usepackage{tikz}
\usepackage{pgfplots}
\usetikzlibrary{shapes.arrows, patterns, calc}
\usepackage{tikz-3dplot}
\usepackage{bbm}
\ifpdf
  \DeclareGraphicsExtensions{.eps,.pdf,.png,.jpg}
\else
  \DeclareGraphicsExtensions{.eps}
\fi

\usepgfplotslibrary{groupplots}
\usetikzlibrary{arrows}
\usetikzlibrary{shapes}
\usetikzlibrary{decorations.text}
\usetikzlibrary{quantikz}
\usepackage{siunitx}
\usetikzlibrary{arrows.meta}

\pgfplotsset{ compat=1.18,
    standard/.style={
    scale only axis,
    width=0.5\textwidth,
    enlarge x limits=0.05,
    enlarge y limits=0.05,
    max space between ticks=40,
    every axis/.append style={font=\normalsize},
	every legend/.append style={font=\normalsize},
	every node/.append style={font=\normalsize},	
	}
}

\definecolor{steelblue}{HTML}{A1BDC7}
\definecolor{orange}{HTML}{D98C21}
\definecolor{silver}{HTML}{B0ABA8}
\definecolor{rust}{HTML}{B8420F}
\definecolor{seagreen}{HTML}{2E6B69}
\definecolor{joshua}{HTML}{FBDC7F}
\definecolor{darksky}{HTML}{154c79}

\colorlet{lightsilver}{silver!30!white}
\colorlet{darkorange}{orange!85!black}
\colorlet{darksilver}{silver!85!black}
\colorlet{darksteelblue}{steelblue!85!black}
\colorlet{darkrust}{rust!85!black}
\colorlet{darkseagreen}{seagreen!85!black}

\usepackage{hyperref}
\usepackage{cleveref}  
\hypersetup{colorlinks=true,linkcolor=darkrust,citecolor=darkseagreen,urlcolor=darksilver}

\usepackage{mathtools}
\usepackage{stackengine}
\usepackage{fixmath}
\usepackage{xcolor}
\usepackage{MnSymbol}
\mathtoolsset{centercolon}  

\renewcommand{\phi}{\varphi}
\renewcommand{\epsilon}{\varepsilon}
\newcommand{\eps}{\varepsilon}

\newcommand{\cnst}[1]{\mathrm{#1}}

\newcommand{\Id}{\mathbf{I}} 

\providecommand{\mathbbm}{\mathbb} 

\newcommand{\R}{\mathbbm{R}}

\newcommand{\abs}[1]{\left|#1\right|}

\newcommand{\vct}[1]{\mathbold{#1}}
\newcommand{\mtx}[1]{\mathbold{#1}}

\newcommand{\ad}{{*}}

\newcommand{\trace}{\operatorname{tr}}

\newcommand{\DExp}[3]{\operatorname{Exp}_{\left(#2\right)}^{#1} \left(#3\right)}

\newcommand{\Exp}[2]{\operatorname{Exp}_{\left(#1\right)}\left(#2\right)}

\newcommand{\Tangent}[2]{T_{#1}{#2}}

\newcommand{\manif}{\mathcal{M}}
\newcommand{\ambient}{\mathcal{A}}
\newcommand{\extended}{\mathcal{E}}

\newcommand{\diff}{\Phi}
\newcommand{\deriv}[1]{#1^\prime}

\newcommand{\bpot}{\psi}
\newcommand{\emb}{\xi}

\DeclarePairedDelimiterX{\infdivx}[2]{(}{)}{%
  #1\;\delimsize\|\;#2%
}
\newcommand{\Div}[1]{\mathbb{D}_{#1}\infdivx}

\newcommand{\divergence}{\operatorname{div}}
\newcommand\D{\mathop{}\cnst{d}}

\newcommand{\defeq}{\vcentcolon=}
\newcommand{\inprod}[1]{\left\langle #1 \right\rangle}

\newcommand{\cdiv}{\overset{\downarrow}{\operatorname{div}}}
\newcommand{\rdiv}{\overset{\rightarrow}{\operatorname{div}}}

\newcommand\Ccancel[2][black]{
    \let\OldcancelColor\CancelColor
    \renewcommand\CancelColor{\color{#1}}
    \cancel{#2}
    \renewcommand\CancelColor{\OldcancelColor}
}

\usepackage{autonum}
\graphicspath{{./figures/tikz}}
\usetikzlibrary{external}
\ifpdf
  \DeclareGraphicsExtensions{.eps,.pdf,.png,.jpg}
\else
  \DeclareGraphicsExtensions{.eps}
\fi

\newsiamremark{remark}{Remark}
\newsiamremark{hypothesis}{Hypothesis}
\crefname{hypothesis}{Hypothesis}{Hypotheses}
\newsiamthm{claim}{Claim}

\headers{Information geometric regularization}{Ruijia Cao and Florian Sch{\"a}fer}

\title{Information geometric regularization \\ of the barotropic Euler equation} 

\author{Ruijia Cao\thanks{Center for Applied Mathematics, Cornell University 
  (\email{rc948@cornell.edu}).}
\and Florian Sch{\"a}fer\thanks{Courant Institute of Mathematical Sciences, New York University
  (\email{florian.schaefer@nyu.edu}).}}

\usepackage{amsopn}

\ifpdf
\hypersetup{
  pdftitle={Information geometric regularization of the barotropic Euler equation},
  pdfauthor={Ruijia Cao, and Florian Sch{\"a}fer}
}
\fi

\begin{document}

\maketitle

\begin{abstract}
  Shock waves in gas dynamics feature jump discontinuities that hinder numerical simulations.  
  Viscous regularizations are prone to excessive dissipation of fine-scale structures.
  In this work, we propose the first inviscid regularization of the multidimensional Euler equation based on ideas from semidefinite programming, information geometry, geometric hydrodynamics, and nonlinear elasticity. 
  The Lagrangian flow maps of Euler solutions are a dynamical system on the manifold of diffeomorphisms. 
  We observe that shock formation arises from the manifold's geodesic incompleteness.
  Our regularization embeds it into an ambient space equipped with the information geometry of the logarithmic barrier function. 
  Thus, the diffeomorphism manifold inherits a geodesically complete geometry.
  The resulting regularized conservation law replaces shocks with smooth profiles without affecting oscillatory structures.
  One and two-dimensional numerical experiments show its practical potential to enable higher-order methods without explicit shock capturing.
  While we focus on the barotropic Euler equations for concreteness and simplicity of exposition, our regularization easily extends to more general Euler and Navier-Stokes-type equations.
  Our approach regularizes the Wasserstein geometry of the mass density with its information geometry.
  The former captures the natural trajectories of physical particles and the latter that of statistical estimators.
  Information geometric regularization accounts for the mass density's dual nature as a statistical/computational tool summarizing the motion of physical particles.
  Thus, our work is a starting point for \emph{information geometric mechanics} that views solutions of continuum mechanical PDEs as parameters of statistical models for unresolved scales and uses their information geometry to evolve them in time.
\end{abstract}



\section{Introduction}
\subsection{The problem} This work is concerned with the treatment of shocks in computational gas dynamics.
Virtually all high-speed gas flows experience a steepening of wavefronts over time, due to faster upstream gas particles catching up to slower downstream ones.
Viscous forces ensure the smoothness wavefronts on microscopic scales of the order of the mean free path of the gas particles \cite{koreeda1995front}. 
But on the scale of any practical computational mesh, wave steepening results in discontinuities of the velocity and density fields, known as shocks. 
This work addresses the fundamental challenge of \emph{shock capturing} by designing a regularized PDE model that coarse-grains the effect of the shock while minimizing dissipative and dispersive errors.

\subsubsection*{The barotropic Euler equation \nopunct}
of gas dynamics serves as our model problem
\vspace{-6mm}
\begin{equation}
     \label{eqn:euler_velocity}
     \partial_t
     \begin{pmatrix}
          \rho \vct{u} \\
          \rho
     \end{pmatrix}
     + \rdiv
     \begin{pmatrix}
          \rho \vct{u} \otimes \vct{u} + P(\rho) \Id \\ 
          \rho \vct{u}
     \end{pmatrix}
     = 
     \begin{pmatrix} 
          \vct{f}\\
          0
     \end{pmatrix}.
     \vspace{-4mm}
\end{equation}
Here, $\rdiv$ denotes the row-wise divergence, $\vct{u}(\vct{x},t)$ the velocity of a gas particle passing through a point $\vct{x}$ at time $t$, and $\rho(\vct{x}, t)$ the gas density in $\vct{x}$ at $t$.
Thus, the momentum density in a point $\vct{x} \in \R^d$ at time $t$ is  $\vct{\mu}(\vct{x}, t) = \rho(\vct{x}, t) \vct{u}(\vct{x}, t)$ .
The pressure $P$ depends only on the density and $\vct{f}$ denotes external forces acting on the gas. 
We choose \eqref{eqn:euler_velocity} for simplicity of exposition, but our approach is not limited to it.

\begin{figure}
       \begin{tikzpicture}
\begin{groupplot}[
	group style={group size=4 by 2,
	horizontal sep=0 cm,
    vertical sep= 0 cm,},
	]
	\nextgroupplot[
		standard,
		no markers,
		height=0.23\textwidth,
		width=0.23\textwidth,
        ymin=-3,
        ymax=3,
        xmin=-0.1,
        xmax=1.1,
		ylabel={$u$},
		ylabel style={rotate=-90},
		ytick={0},
		yticklabels={0},
		ymajorticks=false,
		ymajorgrids=true,
		xmajorticks=false,
		xticklabels={,,},
		enlarge x limits=0.,
		enlarge y limits=0.,
		legend style={
			at={(2.0,0.975)}, inner sep=3pt,anchor=south,legend columns=4,legend cell align={left}, draw=none,fill=none},
		]	

		\addlegendimage{steelblue, ultra thick,}
		\addlegendimage{joshua,ultra thick, dashed}
		\addlegendimage{orange,ultra thick, loosely dotted}

    \legend{Low Viscosity, Medium Viscosity, High Viscosity};

    \addplot[steelblue, ultra thick, smooth] table[x index={1},y index={0}] {csv/vanishing_viscosity_visualization/euler_low_vis_u_0.0.csv};
    \addplot[joshua, ultra thick, smooth, dashed] table[x index={1},y index={0}] {csv/vanishing_viscosity_visualization/euler_mid_vis_u_0.0.csv};
    \addplot[orange, ultra thick, smooth, loosely dotted] table[x index={1},y index={0}] {csv/vanishing_viscosity_visualization/euler_high_vis_u_0.0.csv};

	\nextgroupplot[
		standard,
		no markers,
		height=0.23\textwidth,
		width=0.23\textwidth,
        ymin=-3,
        ymax=3,
        xmin=-0.1,
        xmax=1.1,
        ytick={0},
		ymajorticks=false,
		ymajorgrids=true,
		xmajorticks=false,
		xticklabels={,,}, 
		enlarge y limits=0.,
		]	

       	\addplot[steelblue, ultra thick, smooth] table[x index={1},y index={0}] {csv/vanishing_viscosity_visualization/euler_low_vis_u_0.064.csv};
        \addplot[joshua, ultra thick, smooth, dashed] table[x index={1},y index={0}] {csv/vanishing_viscosity_visualization/euler_mid_vis_u_0.064.csv};
        \addplot[orange, ultra thick, smooth, loosely dotted] table[x index={1},y index={0}] {csv/vanishing_viscosity_visualization/euler_high_vis_u_0.064.csv};

	\nextgroupplot[
		standard,
		height=0.23\textwidth,
		width=0.23\textwidth,
        ymin=-3,
        ymax=3,
        xmin=-0.1,
        xmax=1.1,
        ytick={0},
		ymajorgrids=true,
		ymajorticks=false,
		xmajorticks=false,
		xticklabels={,,},
		enlarge x limits=0.,
		enlarge y limits=0.,
		]	

       	\addplot[steelblue, ultra thick, smooth] table[x index={1},y index={0}] {csv/vanishing_viscosity_visualization/euler_low_vis_u_0.16.csv};
        \addplot[joshua, ultra thick, smooth, dashed] table[x index={1},y index={0}] {csv/vanishing_viscosity_visualization/euler_mid_vis_u_0.16.csv};
        \addplot[orange, ultra thick, smooth, loosely dotted] table[x index={1},y index={0}] {csv/vanishing_viscosity_visualization/euler_high_vis_u_0.16.csv};

	\nextgroupplot[
		standard,
		height=0.23\textwidth,
		width=0.23\textwidth,
        ymin=-3,
        ymax=3.0,
        xmin=-0.1,
        xmax=1.1,
		yticklabel pos = right,
		ytick={0},
		ytick style = {draw= {none}},
		ymajorgrids=true,
		xmajorticks=false,
		xticklabels={,,},
		enlarge x limits=0.,
		enlarge y limits=0.,
		]	

       	\addplot[steelblue, ultra thick, smooth] table[x index={1},y index={0}] {csv/vanishing_viscosity_visualization/euler_low_vis_u_0.32.csv};
        \addplot[joshua, ultra thick, smooth, dashed] table[x index={1},y index={0}] {csv/vanishing_viscosity_visualization/euler_mid_vis_u_0.32.csv};
        \addplot[orange, ultra thick, smooth, loosely dotted] table[x index={1},y index={0}] {csv/vanishing_viscosity_visualization/euler_high_vis_u_0.32.csv};

\end{groupplot}

\node (title) at ($(group c2r1.center)!0.5!(group c3r1.center)+(0,-1.9cm)$) {\Large $\xlongrightarrow[]{\quad \mathrm{time} \quad}$};

\end{tikzpicture}
      \vspace{-4mm}
      \caption{\textbf{Vanishing viscosity solutions.} Solutions with shocks are defined by vanishing viscosity limits. For finite viscosity, the solution smoothes out over time. }
     \label{fig:vanishing_viscosity}
     \vspace{-7mm}
\end{figure}

\subsubsection*{Vanishing viscosity solutions and entropy conditions}
Shock formation in Euler's equations was proved for a wide range of flows \cite{sideris1985formation,christodoulou1993global,christodoulou2007formation,christodoulou2014compressible,christodoulou2019shock}. 
Shock waves preclude long-time existence of classical solutions to \cref{eqn:euler_velocity}. 
The canonical solution concept are limits of solutions with vanishing viscous regularization \cite{guermond2014viscous}(see \cref{fig:vanishing_viscosity}). 
They are also characterized by entropy conditions (see, e.g. \cite{leveque1992numerical}).  

\subsection{Existing approaches}
\label{sec:existing_approaches}
\subsubsection*{Localized viscosity}
Naive approaches for computing vanishing viscosity solutions require a mesh width proportional to the viscosity, causing excessive computational cost or energy dissipation. 
Instead, numerous methods add viscosity \emph{adaptively}, governed by so-called shock sensors \cite{vonneumann1950method,puppo2004numerical,cook2005hyperviscosity,fiorina2007artificial,mani2009suitability,barter2010shock,guermond2011entropy,bruno2022fc,dolejvsi2003some}. 

\subsubsection*{Purely numerical remedies \nopunct}
emphasize the avoidance of Gibbs oscillations encountered by high-order schemes near shocks using limiters that reduce the approximation order at the shock \cite{van1979towards,ray2018artificial,liu1994weighted,harten1997uniformly,shu1998advanced}.
The numerical dissipation of low-order methods result in behavior similar to localized viscosity.
However, a downside of these methods is the absence of a well-defined PDE.
For instance, stencil-switching and slope limiters introduce errors in the derivatives of the solution with respect to problem parameters \cite{bodony2022adjoint,lozano2019watch}.

\subsubsection*{Inviscid PDE-based regularization}
A third line of work searches for \emph{invsicid} PDE-based regularizations.
First attempts based on \cite{leray1934essai} failed to capture the correct shock speeds of the Euler equations, making them unsuitable for most applications \cite{bhat2006hamiltonian,bhat2009riemann,bhat2009regularization}.
The inviscid regularizations of \cite{clamond2018non,guelmame2022hamiltonian} maintain the correct shock speed but are limited to unidimensional problems and can still form singularities \cite{pu2018weakly,guelmame2022global}.

\subsection{This work\nopunct} presents the first inviscid regularization of Euler-like equations.

\subsubsection*{Interior point methods for Euler} 
For vanishing $P, \vct{f}$, \cref{eqn:euler_velocity} is solvable by tracking characteristic curves --- trajectories of fluid particles.
Shocks form when particles collide and merge inelastically.
Equivalently, particle trajectories are minimizers of a variational problem with a non-crossing constraint.
This work modifies characteristics to follow solution paths of interior point methods defined by a barrier $\psi$ \cite{nemirovski2008interior}.
This ensures strict feasibility (shock avoidance) while preserving long-time behavior. 

\subsubsection*{Information geometry on diffeomorphism manifolds}
To make this idea precise, we use the fact that solution paths of interior point methods are dual geodesics in the information geometry defined by $\psi$ \cite{amari2010information,amari2016information}. 
Collectively, characteristics describe geodesics on diffeomorphism manifolds \cite{arnold1966geometrie,ebin1970groups,khesin2007shock,khesin2021geometric}.
Our regularization replaces them with dual geodesics in the information geometry of $\bpot$ (see \cref{eqn:barrier_unidimensional,eqn:extended_barrier} for the definition of $\bpot$).
Adding forces due to $\vct{f}, P$ and returning to a Eulerian description yields the regularized Euler equation (defining the Jacobian matrix $[\cnst{D}\vct{u}]_{ij} \coloneqq \partial_{j} u_i$),
\begin{equation}
     \label{eqn:reg_euler_intro}
     \begin{cases}
          \partial_t
          \begin{pmatrix}
               \rho \vct{u} \\
               \rho
          \end{pmatrix}
          + \rdiv 
          \begin{pmatrix}
               \rho \vct{u} \otimes \vct{u} + \left(P(\rho) + \Sigma\right) \Id\\ 
               \rho \vct{u}
          \end{pmatrix}
          = 
          \begin{pmatrix} 
               \vct{f}\\
               0
          \end{pmatrix} \\
          \rho^{-1} \Sigma - \alpha \divergence(\rho^{-1} \nabla \Sigma) = \alpha  \left(\trace^2\left([\cnst{D} \vct{u}]\right) + \trace\left([\cnst{D} \vct{u}]^2\right) \right).
     \end{cases}
     \vspace{-1mm}
\end{equation}
This work focuses on the barotropic Euler equations for concreteness and to simplify the exposition. 
But one can easily derive extensions to general Euler and Navier-Stokes-type equations by adding the \emph{entropic pressure} $\Sigma$ to the mechanical pressure $P$.
This has recently enabled compressible flow simulations at unprecedented scale \cite{wilfong2025simulating}.
By choosing the parameter $\sqrt{\alpha}$ proportional to the mesh size, the elliptic problem defining $\Sigma$ in \cref{eqn:reg_euler_intro} is solvable by a constant number of Jacobi or Gauss-Seidel sweeps.
Two sweeps per flux computation suffice in practice, at negligible computational cost.

\section{Background: Four perspectives on shock formation}
\label{sec:background}
\subsection{Model problem: The unidimensional pressureless case}
The unidimensional pressureless case of \cref{eqn:euler_velocity}, with $P, \vct{f} \equiv 0$ and domain $\R$, is
\vspace{-1mm}
\begin{equation}
     \label{eqn:burger_velocity}
     \partial_t
     \begin{pmatrix}
          \rho u \\
          \rho
     \end{pmatrix}
     + \partial_x
     \begin{pmatrix}
          \rho u^2 \\ 
          \rho u
     \end{pmatrix}
     = 
     \begin{pmatrix} 
          0\\
          0
     \end{pmatrix}, \text{with initial conditions $\rho(\cdot, 0) \equiv 1$ and $u(\cdot, 0) \equiv u_0(\cdot)$.}
     \vspace{-1mm}
\end{equation}
Smooth solutions of \cref{eqn:burger_velocity} describe the inertial movement of noninteracting particles. 
We will now discuss four closely related perspectives on shock formation in this system as summarized in \cref{tab:perspectives}.
\begin{table}
     \centering
     \begin{tabular}{m{0.1335\textwidth}|m{0.167\textwidth}|m{0.195\textwidth}|m{0.1675\textwidth}|m{0.167\textwidth}|}
                      \hline
                       Perspective& Trajectories, \Cref{sec:shock_trajectories} & Diffeomorphism, \Cref{sec:shock_diffeomorphisms} & Transport, \newline\Cref{sec:shock_transport} & Optimization, \Cref{sec:shock_optimization} \\
                      \hline
                      Shock \newline formation & Particles\newline collide. & Flow maps reach the manifold's \newline boundary.& Pushforward measure forms atoms. & Monotonicity constraint activates. \\ 
                      \hline
                      Entropy condition & Particles stick to each other & Flow maps stick to the manifold's boundary& Formation of atoms is never undone. & First-order/KKT conditions. \\
                      \hline
     \end{tabular}
     \vspace{-2mm}
     \caption{Four perspectives on unidimensional pressureless shock formation.}
     \vspace{-8mm}
     \label{tab:perspectives}
\end{table}
We solve \cref{eqn:burger_velocity} by tracking particle trajectories. 
This approach coincides with the method of characteristics \cite{evans2022partial}. 
For initial conditions $u(\cdot, 0)$ we define a family of curves indexed by $x \in \R$ as $\left\{t \mapsto \diff_t(x) \coloneqq x + t u_0(x)\right\}_{x \in \R}$.
Without particle interactions, $t \mapsto \diff_{t}(x)$ is the position at time $t$ of the particle initially located at $x$, and $\dot{\diff}_t(x)$ is its derivative. 
We obtain the velocity $u(y, t)$ in position $x$ at time $t$ by tracing the particle passing through $y$ at time $t$ back to its initial position $x$ at time $0$.
Likewise, we determine the density $\rho(y, t)$ by measuring whether nearby particles tend to disperse or accumulate, which is reflected by the spatial derivative $\partial_x \diff_t(x)$. 
\begin{minipage}{0.3825\textwidth}
     This yields the explicit solution\newline
     \vspace{-2mm}
     \begin{equation}
          \label{eqn:characteristic_solution}
     \begin{pmatrix}
          {\color{rust} u(\diff_{t}(x), t)} \\
          {\color{orange}\rho(\diff_{t}(x), t)}
        \end{pmatrix}
        \defeq
        \begin{pmatrix}
          {\color{rust} u_0(x)}\\
          {\color{orange}(\partial_x \diff_{t}(x))^{-1}}
        \end{pmatrix}.
     \end{equation}
     \vspace{-5mm}
     \captionof{figure}{\textbf{Solution by particle-tracing.} Tracing trajectories to the initial condition provides $u$. Computing their dispersion or convervence yields $\rho$.}
\end{minipage}
\begin{minipage}{0.610\textwidth}
        \begin{tikzpicture}
	\begin{scope}[xshift = -0.00\textwidth]
		\draw[silver, ultra thick, ->] (-0.5, 0) -- (3, 0);
		\draw[silver, ultra thick, ->] (0, -0.5) -- (0, 3);



		\draw[rust, ultra thick, ->] (1, 0) -- (2, 2.5);

		\node[] at (1,-0.375) {$x$};
		\draw[silver, very thick] (1,-0.125) -- (1,0.125);
		\node[] at (-0.25,2.5) {$t$};
		\node[] at (2.725,2.725) {$(\diff_t(x), t)$};
		\draw[silver, very thick] (-0.125,2.5) -- (0.125,2.5);
		\node[] at (2.075,-0.375) {$\diff_t(x)$};
		\draw[silver, very thick] (2,-0.125) -- (2,0.125);
	\end{scope}

	\node[draw=steelblue, ultra thick] at (3.650,-1.00) {$\diff_t(x) \coloneqq x + t u_0(x)$};

	\begin{scope}[xshift = 0.480\textwidth]
		\draw[silver, ultra thick, ->] (-0.5, 0) -- (3, 0);
		\draw[silver, ultra thick, ->] (0, -0.5) -- (0, 3);

		\draw[orange, dashed, ultra thick, ->] (0.80, 0) -- (1.95, 2.5);
		\draw[orange, dashed, ultra thick, ->] (0.60, 0) -- (1.85, 2.5);
		\draw[orange, dashed, ultra thick, ->] (0.40, 0) -- (1.75, 2.5);
		\draw[orange, dashed, ultra thick, ->] (0.20, 0) -- (1.65, 2.5);

		\draw[orange, dashed, ultra thick, ->] (1.2, 0) -- (2.2, 2.5);
		\draw[orange, dashed, ultra thick, ->] (1.4, 0) -- (2.5, 2.5);
		\draw[orange, dashed, ultra thick, ->] (1.6, 0) -- (2.9, 2.5);
		\draw[orange, dashed, ultra thick, ->] (1.8, 0) -- (3.5, 2.5);

		\draw[silver, ultra thick, ->] (1, 0) -- (2, 2.5);

		\node[] at (1,-0.375) {$x$};
		\draw[silver, very thick] (1,-0.125) -- (1,0.125);
		\node[] at (-0.25,2.5) {$t$};
		\node[] at (2.725,2.725) {$(\diff_t(x), t)$};
		\draw[silver, very thick] (-0.125,2.5) -- (0.125,2.5);
		\node[] at (2.075,-0.375) {$\diff_t(x)$};
		\draw[silver, very thick] (2,-0.125) -- (2,0.125);
	\end{scope}

\end{tikzpicture}
\end{minipage}
\newline

\begin{figure}[h]
     \begin{tikzpicture}
\begin{groupplot}[
	group style={group size=4 by 1,
	horizontal sep=0 cm,
    vertical sep= 0 cm,},
	]
	\nextgroupplot[
		standard,
		no markers,
		height=0.230\textwidth,
		width=0.230\textwidth,
        ymin=-0.8,
        ymax=1.4,
        xmin=0.2,
        xmax=0.8,
		ylabel={{\color{seagreen} $u$}, {\color{steelblue} $t$}},
		xlabel={{\color{seagreen} $x$}, {\color{steelblue} $\diff_t(x)$}},
		ytick={0},
		yticklabels={0},
		ymajorticks=false,
		ymajorgrids=true,
		xmajorticks=false,
		xticklabels={,,},
		enlarge x limits=0.,
		enlarge y limits=0.,
		legend style={
			at={(2.0,1.075)}, inner sep=3pt,anchor=south,legend columns=4,legend cell align={left}, draw=none,fill=none},
		]	

		\addlegendimage{blue,very thick}
		\addlegendimage{orange,very thick}
		\addlegendimage{green,very thick}

        \foreach \col in {1,...,25}{
            \addplot[steelblue, very thick, smooth] table[x index = {\col}, y expr={6 * (\thisrowno{0})}] {csv/burgers_characteristics/burgers_chars_0.01.csv};
        }

    \addplot[seagreen, ultra thick, smooth] table[x index={1},y index={0}] {csv/burgers_characteristics/burgers_u_0.01.csv};

	\nextgroupplot[
		standard,
		no markers,
		height=0.230\textwidth,
		width=0.230\textwidth,
        ymin=-0.8,
        ymax=1.4,
        xmin=0.2,
        xmax=0.8,
        ytick={0},
		ymajorticks=false,
		ymajorgrids=true,
		xmajorticks=false,
		xticklabels={,,}, 
		enlarge y limits=0.,
		legend style={
			at={(2.20,1.00)},inner sep=3pt,anchor=south,legend columns=4,legend cell align={left}, draw=none,fill=none},
		]	

        \foreach \col in {1,...,25}{
            \addplot[steelblue, very thick, smooth] table[x index = {\col}, y expr={6 * (\thisrowno{0})}] {csv/burgers_characteristics/burgers_chars_0.05.csv};
        }

       	\addplot[seagreen, ultra thick, smooth] table[x index={1},y index={0}] {csv/burgers_characteristics/burgers_u_0.05.csv};

	\nextgroupplot[
		standard,
		height=0.230\textwidth,
		width=0.230\textwidth,
        ymin=-0.8,
        ymax=1.4,
        xmin=0.2,
        xmax=0.8,
        ytick={0},
		ymajorgrids=true,
		ymajorticks=false,
		xmajorticks=false,
		xticklabels={,,},
		enlarge x limits=0.,
		enlarge y limits=0.,
		legend style={
			at={(2.20,1.00)},inner sep=3pt,anchor=south,legend columns=4,legend cell align={left}, draw=none,fill=none},
		]	

        \foreach \col in {1,...,25}{
            \addplot[steelblue, very thick, smooth] table[x index = {\col}, y expr={6 * (\thisrowno{0})}] {csv/burgers_characteristics/burgers_chars_0.1.csv};
        }

       	\addplot[seagreen, ultra thick, smooth] table[x index={1},y index={0}] {csv/burgers_characteristics/burgers_u_0.1.csv};

	\nextgroupplot[
		standard,
		height=0.230\textwidth,
		width=0.230\textwidth,
        ymin=-0.8,
        ymax=1.4,
        xmin=0.2,
        xmax=0.8,
		yticklabel pos = right,
		ytick={0},
		ytick style = {draw= {none}},
		ymajorgrids=true,
		xmajorticks=false,
		xticklabels={,,},
		enlarge x limits=0.,
		enlarge y limits=0.,
		legend style={
			at={(2.20,1.00)},inner sep=3pt,anchor=south,legend columns=4,legend cell align={left}, draw=none,fill=none},
		]	

        \foreach \col in {1,...,25}{
            \addplot[orange, very thick, smooth] table[x index = {\col}, y expr={6 * (\thisrowno{0})}] {csv/burgers_characteristics/burgers_chars_0.159.csv};
        }
       	\addplot[rust, ultra thick, smooth,] table[x index={1},y index={0}] {csv/burgers_characteristics/burgers_u_0.159.csv};

\end{groupplot}

\node (title) at ($(group c2r1.center)!0.5!(group c3r1.center)+(0,-1.9cm)$) {\Large $\xlongrightarrow[]{\quad \mathrm{time} \quad}$};

\end{tikzpicture}
     \vspace{-6mm}
     \caption{\textbf{Perspective I: Shock formation and characteristics.} The plot overlays particle trajectories $x \mapsto \diff_t(x)$ with the velocity of a solution of \cref{eqn:burger_velocity}, at four different time steps. As the trajectories approach each other, the profile steepens. In the rightmost plot, the trajectories cross and the flow develops a shock. }
     \label{fig:characteristic_shock}
     \vspace{-2mm}
\end{figure}

\subsection{Perspective I: Particle trajectories}
\label{sec:shock_trajectories}
So long as the particle trajectories do not cross, the closed form solution \cref{eqn:characteristic_solution} is well-defined, as there is at most one particle passing through each point. 
But for the vast majority of initial conditions, there exists a time $T$ at which two trajectories cross, meaning that for $x_1 \neq x_2$, we have $\diff_T(x_1) = \diff_T(x_2)$. 
Since the two particles have different velocities, the values of $u$ to the left and right of $\diff_T(x_1) = \diff_T(x_2)$ differ. 
A shock discontinuity has formed, as illustrated in \cref{fig:characteristic_shock}.
Starting at $T$, \cref{eqn:burger_velocity} ceases to have classical solutions.
Considering vanishing viscosity solutions amounts to assuming that particles are ``sticky'' and merge upon collision, as shown in \cref{fig:collision}. 
The velocity of the combined particles is chosen such as to conserve linear momentum. E.g. \cite{cavalletti2015simple,dermoune1999probabilistic,brenier2013sticky,carrillo2023equivalence,suder2021lagrangian,hynd2019lagrangian,hynd2020trajectory,nguyen2008pressureless} study this system.
\begin{figure}[h]
     \begin{tikzpicture}
\begin{groupplot}[
	group style={group size=4 by 1,
	horizontal sep=0.60cm,
    vertical sep=0.5cm,},
	]
	\nextgroupplot[
		standard,
 		xlabel={$\diff_t(x)$},
		ylabel style={yshift=-6.0pt},	
		ylabel={$t$},
		height=0.2025\textwidth,
		width=0.2025\textwidth,
        ymin=-0.05,
        ymax=1.0,
        xmin=-0.05,
        xmax=1.0,
        xticklabels={{$x_1$}, {$x_2$}},
        yticklabels={{$0$}},
		xtick={0,0.5},
		ytick={0},
		enlarge x limits=0.,
		enlarge y limits=0.,
		]	

        \addplot[steelblue,very thick] table[x index={1},y index={0}, col sep=comma] {csv/regularization/characteristics_alpha_0.0_medium.txt};
        \addplot[steelblue,very thick] table[x index={2},y index={0}, col sep=comma] {csv/regularization/characteristics_alpha_0.0_medium.txt};
		
		\draw[darksilver] (axis cs: 0.0,0.0) -- (axis cs: 1.5,0.0);

	\nextgroupplot[
		standard,
 		xlabel={$\diff_t(x)$},
		ylabel={},
		height=0.2025\textwidth,
		width=0.2025\textwidth,
        ymin=-0.05,
        ymax=1.0,
        xmin=-0.05,
        xmax=1.0,
        xticklabels={{$x_1$}, {$x_2$}},
        yticklabels={},
		xtick={0,0.5},
		ytick={0},
		enlarge x limits=0.,
		enlarge y limits=0.,
		]	

        \addplot[orange,very thick] table[x index={1},y index={0}, col sep=comma] {csv/regularization/characteristics_alpha_0.0_long.txt};
        \addplot[orange,very thick] table[x index={2},y index={0}, col sep=comma] {csv/regularization/characteristics_alpha_0.0_long.txt};
		
		\draw[darksilver] (axis cs: 0.0,0.0) -- (axis cs: 1.5,0.0);

\nextgroupplot[
		standard,
 		xlabel={$\diff_t(x_1)$},
		ylabel={$\diff_t(x_2)$},
		height=0.2025\textwidth,
		width=0.2025\textwidth,
        ymin=0.775,
        ymax=2.475,
        xmin=0.775,
        xmax=2.475,
		xtick={0.775},
        xticklabels={\phantom{{$x_1$}, {$x_2$}}},
        yticklabels={},
		ymajorticks=false,
		enlarge x limits=0.,
		enlarge y limits=0.,
		]	

		\draw[darksilver] (axis cs: -0.5,-0.5) -- (axis cs: 2.5,2.5);

       	\addplot[steelblue,very thick] table[x expr={1 + \thisrowno{1} - \thisrowno{2}},y expr={ 1 + \thisrowno{1} + \thisrowno{2}}, col sep=comma] {csv/regularization/phase_space_alpha_0.0_medium.txt};

		\addplot+[mark=none,
		draw=none,
        domain=-0.0:2.50,
        samples=100,
        pattern= dots,
        area legend,
        pattern color=darksilver]{x} \closedcycle;

\nextgroupplot[
		standard,
 		xlabel={$\diff_t(x_1)$},
		ylabel={},
		height=0.2025\textwidth,
		width=0.2025\textwidth,
        ymin=0.775,
        ymax=2.475,
        xmin=0.775,
        xmax=2.475,
		xtick={0.775},
        xticklabels={\phantom{{$x_1$}, {$x_2$}}},
        yticklabels={},
		ymajorticks=false,
		enlarge x limits=0.,
		enlarge y limits=0.,
		]	

		\draw[darksilver] (axis cs: -0.5,-0.5) -- (axis cs: 2.5,2.5);

       	\addplot[orange,very thick] table[x expr={1 + \thisrowno{1} - \thisrowno{2}},y expr={ 1 + \thisrowno{1} + \thisrowno{2}}, col sep=comma] {csv/regularization/phase_space_alpha_0.0_long.txt};

		\addplot+[mark=none,
		draw=none,
        domain=-0.0:2.50,
        samples=100,
        pattern= dots,
        area legend,
        pattern color=darksilver]{x} \closedcycle;
      
	\end{groupplot}
\end{tikzpicture}
     \vspace{-6mm}
     \caption{\textbf{Perspective II: The boundary of diffeomorphism manifolds.} Instead of tracking particle trajectories, we can view $t \mapsto \diff_t$ as a geodesic on a diffeomorphism manifold.
     Shock form when the geodesic reaches the manifold's boundary. Vanishing viscosity solutions continue to evolve on the boundary. }
     \label{fig:collision}
     \vspace{-2mm}
\end{figure}
\subsection{Perspective II: Diffeomorphism manifolds}
\label{sec:shock_diffeomorphisms}
Instead of the trajectories $t \mapsto \diff_t(x)$ of individual particles, we can consider the function-valued trajectory $t \mapsto \diff_t$ of the so-called flow maps.
In \cite{arnold1966geometrie}, Arnold observed that flow maps of solutions to the incompressible Euler solutions are geodesic curves on the manifold of volume-preserving diffeomorphisms, using the $L^2$ inner product of vector fields as a Riemannian metric. 
In the compressible Euler equations, the flow maps instead evolve on the manifold of \emph{general} diffeomorphisms \cite{khesin2021geometric}. 
Smooth solutions of pressureless Euler equations are $L^2$-geodesics of the form $\diff_t = \diff_0 + t \vct{v}$ for a vector field $\vct{v}$.
The $L^2$-geometry of the manifold of general diffeomorphisms is not geodesically complete, meaning that geodesic curves do not continue indefinitely.
When two particle trajectories cross, $\diff_t$ loses injectivity and ceases to be a diffeomorphism.
By tracking, as in \cref{fig:collision}, the trajectories of just two particles, we can visualize two-dimensional projections of the diffeomorphism manifold.
The merging of particles in vanishing viscosity solutions corresponds to the flow map moving along the boundary of the diffeomorphism manifold.

\begin{figure}[h]
     \centering
     \begin{tikzpicture}[scale=0.935]

\begin{scope}[xshift=0.001\textwidth]
\def\radius{1}
\def\height{2}
\def\circheight{1.25}

\fill[lightsilver, opacity=0.4] 
  (\radius, \height) -- 
  (\radius, \height) arc[start angle=-0, end angle=-180, x radius=\radius cm, y radius=0.5cm] --
  (\radius, \height) arc[start angle=0, end angle=180, x radius=\radius cm, y radius=0.5cm];

\draw[thick, rust] (\radius, \circheight) arc[start angle=0, end angle=180, x radius=\radius cm, y radius=0.5cm];

\draw[thick, darksky] (\radius, 0) arc[start angle=0, end angle=180, x radius=\radius cm, y radius=0.5cm];

\fill[steelblue, opacity=0.6] 
  (\radius, \height) -- 
  (\radius, 0) arc[start angle=-0, end angle=-180, x radius=\radius cm, y radius=0.5cm] --
  (-\radius, \height) arc[start angle=-180, end angle=-0, x radius=\radius cm, y radius=0.5cm];

\draw[thick, darksky] (\radius, 0) arc[start angle=-0, end angle=-180, x radius=\radius cm, y radius=0.5cm];

\draw[thick, darksky] (\radius, 0) -- (\radius, \height);
\draw[thick, darksky] (-\radius, 0) -- (-\radius, \height);

\draw[thick, darksky] (0, \height) ellipse (\radius cm and 0.5cm);

\draw[thick, rust] (\radius, \circheight) arc[start angle=-0, end angle=-180, x radius=\radius cm, y radius=0.5cm];
\end{scope}

\begin{scope}
\def\circheight{1.25}
\node at (0.110\textwidth, \circheight) {\color{rust} \Huge $\rightarrow$};
\end{scope}

\begin{scope}[xshift = 0.220\textwidth]
\def\radius{1}
\def\height{2}
\def\circheight{1.25}
\draw[thick, rust] (\radius, \circheight) arc[start angle=0, end angle=180, x radius=\radius cm, y radius=0.5cm];
\draw[thick, rust] (\radius, \circheight) arc[start angle=-0, end angle=-180, x radius=\radius cm, y radius=0.5cm];
\end{scope}

\begin{scope}[xshift = 0.400\textwidth]
	\begin{scope}
		\fill[steelblue, opacity=0.6] (-1.2, 0.5) -- (2.8, 0.5) -- (2.8, 3) -- (-1.2, 3) -- cycle;

		\draw[thick, darksky] (-1.2, 0.5) -- (2.8, 0.5);

		\draw[->, thick] (-0.5, 0) -- (2.5, 0);
		\draw[->, thick] (0, -0.5) -- (0, 2.5);

		\draw[] (-1.2, -1.0) -- (2.8, -1.0) -- (2.8, 3) -- (-1.2, 3) -- cycle;
		\node at (-0.2,0.30) {\color{darksky} $0$};

		\draw[thick, rust] (-1.2, 1.25) -- (2.8, 1.25);
		\node at (-0.6,0.95) {\color{rust} $x_2 - x_1$};

		\node at (1.6,- 0.4) {$\diff(x_2) + \diff(x_1)$};
		\node at (-0.0, 2.7) { $\diff(x_2) - \diff(x_1)$};
	\end{scope}

	\begin{scope}
	\def\circheight{1.25}
	\node at (0.245\textwidth, \circheight) {\color{rust} \Huge $\rightarrow$};
	\end{scope}

	\begin{scope}[xshift = 0.370\textwidth]
		\draw[->, thick] (-0.5, 0) -- (2.5, 0);
		\draw[->, thick] (0, -0.5) -- (0, 2.5);

		\draw[] (-1.2, -1.0) -- (2.8, -1.0) -- (2.8, 3) -- (-1.2, 3) -- cycle;

		\draw[thick, rust] (-1.2, 1.25) -- (2.8, 1.25);
		\node at (-0.6,0.95) {\color{rust} $x_2 - x_1$};


		\node at (1.6,- 0.4) {$\diff(x_2) + \diff(x_1)$};
		\node at (-0.0, 2.7) { $\diff(x_2) - \diff(x_1)$};
	\end{scope}
\end{scope}

\end{tikzpicture}
     \vspace{-2mm}
     \caption{\textbf{Geodesically complete submanifolds.} The lateral surface of a cylinder is not geodesically complete. 
     But the circle, a submanifold of it, is. 
     Similarly, the manifold of unidimensional diffeomorphisms (where $\diff(x_2) - \diff(x_1) > 0$) is not geodesically complete. 
     But its volume-preserving submanifold (where $\diff(x_2) - \diff(x_1) = x_2 - x_1$) is.}
     \label{fig:submanifold_completion}
\end{figure}

\begin{remark}
    Geodesic completeness of volume preserving diffeomorphism manifolds implies the existence of global solutions of the incompressible Euler equations \cite{ebin1970groups}.
    The recent breakthrough of \cite{chen2022stable,chen2023stable} provides a counterexample on bounded domains.
    But this singularity formation is unrelated to shocks and lacks a geometric interpretation like \cref{fig:collision}.
    The difficulty of proving these results suggests that the question of geodesic (in)completeness of manifolds of volume-preserving diffeomorphisms is much more subtle than that of general diffeomorphism manifolds. 
    \cref{fig:submanifold_completion} shows the geometric principle by which incompressibility prevents shock formation.
\end{remark}

\begin{figure}[h]
     \begin{tikzpicture}
\begin{groupplot}[
	group style={group size=4 by 1,
	horizontal sep=0 cm,
    vertical sep= 0 cm,},
	]
	\nextgroupplot[
		standard,
		no markers,
		height=0.230\textwidth,
		width=0.230\textwidth,
        ymin=-0.8,
        ymax=1.4,
        xmin=0.2,
        xmax=0.8,
		ylabel={{\color{seagreen} $\log(\rho)$}, {\color{steelblue} $t$}},
		xlabel={{\color{seagreen} $x$}, {\color{steelblue} $\diff_t(x)$}},
		ytick={0},
		yticklabels={0},
		ymajorticks=false,
		ymajorgrids=true,
		xmajorticks=false,
		xticklabels={,,},
		enlarge x limits=0.,
		enlarge y limits=0.,
		legend style={
			at={(2.0,1.075)}, inner sep=3pt,anchor=south,legend columns=4,legend cell align={left}, draw=none,fill=none},
		]	

		\addlegendimage{blue,very thick}
		\addlegendimage{orange,very thick}
		\addlegendimage{green,very thick}

        \foreach \col in {1,...,25}{
            \addplot[steelblue, very thick, smooth] table[x index = {\col}, y expr={6 * (\thisrowno{0})}] {csv/burgers_characteristics/burgers_chars_0.01.csv};
        }

    \addplot[seagreen, ultra thick, smooth] table[x index={0},y index={1}, y expr={log10(\thisrowno{1})}] {csv/burgers_characteristics/burgers_rho_0.01.csv};

	\nextgroupplot[
		standard,
		no markers,
		height=0.230\textwidth,
		width=0.230\textwidth,
        ymin=-0.8,
        ymax=1.4,
        xmin=0.2,
        xmax=0.8,
        ytick={0},
		ymajorticks=false,
		ymajorgrids=true,
		xmajorticks=false,
		xticklabels={,,}, 
		enlarge y limits=0.,
		legend style={
			at={(2.20,1.00)},inner sep=3pt,anchor=south,legend columns=4,legend cell align={left}, draw=none,fill=none},
		]	

        \foreach \col in {1,...,25}{
            \addplot[steelblue, very thick, smooth] table[x index = {\col}, y expr={6 * (\thisrowno{0})}] {csv/burgers_characteristics/burgers_chars_0.05.csv};
        }

       	\addplot[seagreen, ultra thick, smooth] table[x index={0},y index={1}, y expr={log10(\thisrowno{1})}] {csv/burgers_characteristics/burgers_rho_0.05.csv};

	\nextgroupplot[
		standard,
		height=0.230\textwidth,
		width=0.230\textwidth,
        ymin=-0.8,
        ymax=1.4,
        xmin=0.2,
        xmax=0.8,
        ytick={0},
		ymajorgrids=true,
		ymajorticks=false,
		xmajorticks=false,
		xticklabels={,,},
		enlarge x limits=0.,
		enlarge y limits=0.,
		legend style={
			at={(2.20,1.00)},inner sep=3pt,anchor=south,legend columns=4,legend cell align={left}, draw=none,fill=none},
		]	

        \foreach \col in {1,...,25}{
            \addplot[steelblue, very thick, smooth] table[x index = {\col}, y expr={6 * (\thisrowno{0})}] {csv/burgers_characteristics/burgers_chars_0.1.csv};
        }

       	\addplot[seagreen, ultra thick, smooth] table[x index={0},y index={1}, y expr={log10(\thisrowno{1})}] {csv/burgers_characteristics/burgers_rho_0.1.csv};

	\nextgroupplot[
		standard,
		height=0.230\textwidth,
		width=0.230\textwidth,
        ymin=-0.8,
        ymax=1.4,
        xmin=0.2,
        xmax=0.8,
		yticklabel pos = right,
		ytick={0},
		yticklabels={0},
		ytick style = {draw= {none}},
		ymajorgrids=true,
		xmajorticks=false,
		xticklabels={,,},
		enlarge x limits=0.,
		enlarge y limits=0.,
		legend style={
			at={(2.20,1.00)},inner sep=3pt,anchor=south,legend columns=4,legend cell align={left}, draw=none,fill=none},
		]	

        \foreach \col in {1,...,25}{
            \addplot[orange, very thick, smooth] table[x index = {\col}, y expr={6 * (\thisrowno{0})}] {csv/burgers_characteristics/burgers_chars_0.159.csv};
        }
       	\addplot[rust, ultra thick, smooth] table[x index={0},y index={1}, y expr={log10(\thisrowno{1})}] {csv/burgers_characteristics/burgers_rho_0.159.csv};
\end{groupplot}

\node (title) at ($(group c2r1.center)!0.5!(group c3r1.center)-(0,1.90cm)$) {\Large $\xlongrightarrow[]{\quad \mathrm{time} \quad}$};

\end{tikzpicture}
     \vspace{-6mm}
     \caption{\textbf{Perspective III: Shock formation and mass transport.} Shock formation in the pressureless Euler equation amounts to the formation of atoms in the pushforward measure. 
     Vanishing viscosity solutions continue to transport the resulting singular distributions.}
     \label{fig:characteristic_shock_rho}
     \vspace{-2mm}
\end{figure}

\subsection{Perspective III: Mass transport}
\label{sec:shock_transport}
By \cref{eqn:characteristic_solution}, the mass density $\rho$ satisfies $\rho(\diff_t(x), t) = (\partial_x \diff_t(x))^{-1}$ before shock formation. 
It is the pushforward of the Lebesgue measure by $\diff_t$, meaning that $\int_{\diff_{t}(x_1)}^{\diff_t(x_2)} \rho(y, t) \D y = \int_{x_1}^{x_2} \D x$ for any $x_1, x_2 \in \R$.
Curves $t \mapsto \diff_t$ on diffeomorphism manifolds define a curve on the manifold of positive measures. 
Pushforwards of $L^2$ geodesics on diffeomorphism manifolds correspond to $L^2$ Wasserstein geodesics on the manifold of positive measures \cite{benamou2000computational,otto2001geometry,ambrosio2008gradient}.
Optimal transport focuses on computing geodesics between endpoints.
Instead, \cref{eqn:burger_velocity} computes $L^2$ (Wasserstein) geodesics starting in $\rho(\cdot, 0)$ in an initial direction prescribed by $u(\cdot, 0)$.
At shock formation there exist $x_1 < x_2$ such that $\diff_T(x_1) = \diff_T(x_2)$, meaning that the pushforward of the Lebesgue measure under $\diff_T$ has a point mass of weight $\int_{x_1}^{x_2} \D x$ at $\diff_T(x_1) = \diff_T(x_2)$, as illustrated in \cref{fig:characteristic_shock_rho}.
Shock formation arises from the geodesic incompleteness of the manifold of absolutely continuous measures. 
The analog of particles sticking together and flow maps $\diff_t$ sticking to the diffeomorphism manifold's boundary is the fact that pushforwards of positive measures create or merge, but never remove atoms.

\begin{figure}
     \pgfplotsset{
    colormap={rusts}{
        color=(lightsilver)
        color=(rust)
    }, 
    colormap={joshuas}{
        color=(lightsilver)
        color=(joshua)
    }, 
    colormap={darkskys}{
        color=(lightsilver)
        color=(darksky)
    }, 
    colormap={seagreens}{
        color=(lightsilver)
        color=(seagreen)
    }, 
}

\begin{tikzpicture}
\begin{groupplot}[
	group style={group size=4 by 1,
	horizontal sep=0.115cm,
    vertical sep=0.5cm,},
	]
	\nextgroupplot[
		standard,
		title={Nominal, $t = 0.05$},
 		xlabel={$\diff_t$},
 		ylabel={$\partial_x \diff_t$},
		height=0.225\textwidth,
		width=0.225\textwidth,
        ymin=-0.425,
        ymax=0.425,
        xmin=0.2,
        xmax=1.05,
		xtick={0.575},
        xticklabels={\phantom{{$x_1$}, {$x_2$}}},
		xmajorticks=false,
		ymajorticks=false,
		zmajorticks=false,
		ytick={0},
		yticklabels={$0\:$},
		enlarge x limits=0.,
		enlarge y limits=0.,
		view={0}{90},
		colormap name=seagreens,
		]	

		\def\t{0.05}
		\addplot3[
    		contour lua,
			contour/levels={0.001, 0.004,0.009,0.016, 0.025, 0.036, 0.049, 0.064, 0.081, 0.1, 0.121, 0.144, 0.169, 0.196, 0.225, 0.256, 0.289, 0.324, 0.361, 0.4, 0.441, 0.484, 0.529, 0.576,},
			contour/labels = false,
			samples=75,
			handler/.style=smooth,
			domain = 0.0:1.1,
			thick,
			y domain = -0.5:0.5,
		]{(\t ^ 2) * (1.5^2 + 1^2)/ 2  - \t  * (3/2 * (x - 0.25) - 1 * (y - 0.25)) + 1 / 2 * ((x - 0.25)^2 + (y - 0.25)^2)};

		\draw[darksilver] (axis cs: 0.0,0.0) -- (axis cs: 1.5,0.0);

       	\addplot[very thick, dashed, lightsilver] table[x index={1},y index={2}, col sep=comma] {csv/regularization/phase_space_alpha_0.0.txt};
       	\addplot[steelblue,very thick] table[x index={1},y index={2}, col sep=comma] {csv/regularization/phase_space_alpha_0.0_short.txt};

		\addplot+[mark=none,
		draw=none,
       domain=-0.0:2.50,
       samples=100,
       pattern= dots,
       area legend,
       pattern color=darksilver]{-10} \closedcycle;

	\nextgroupplot[
		standard,
		title={Nominal, $t = 0.45$},
 		xlabel={$\diff_t$},
		height=0.225\textwidth,
		width=0.225\textwidth,
        ymin=-0.425,
        ymax=0.425,
        xmin=0.2,
        xmax=1.05,
		ylabel shift=-125pt,
		xtick={0.575},
        xticklabels={\phantom{{$x_1$}, {$x_2$}}},
		yticklabel pos=right,
		xmajorticks=false,
		ymajorticks=false,
		zmajorticks=false,
		ytick={0},
		yticklabels={$0\:$},
		enlarge x limits=0.,
		enlarge y limits=0.,
		view={0}{90},
		colormap name=rusts,
		]	

		\def\t{0.45}
		\addplot3[
    		contour lua,
			contour/levels={0.001, 0.004,0.009,0.016, 0.025, 0.036, 0.049, 0.064, 0.081, 0.1, 0.121, 0.144, 0.169, 0.196, 0.225, 0.256, 0.289, 0.324, 0.361, 0.4, 0.441, 0.484, 0.529, 0.576, 0.625, 0.676, 0.729, 0.784,},
			contour/labels = false,
			samples=75,
			handler/.style=smooth,
			domain = 0.0:1.1,
			thick,
			y domain = -0.5:0.5,
		]{(\t ^ 2) * (1.5^2 + 1^2)/ 2  - \t  * (3/2 * (x - 0.25) - 1 * (y - 0.25)) + 1 / 2 * ((x - 0.25)^2 + (y - 0.25)^2)};

		\draw[darksilver] (axis cs: 0.0,0.0) -- (axis cs: 1.5,0.0);

       	\addplot[very thick, dashed, lightsilver] table[x index={1},y index={2}, col sep=comma] {csv/regularization/phase_space_alpha_0.0.txt};
       	\addplot[orange,very thick] table[x index={1},y index={2}, col sep=comma] {csv/regularization/phase_space_alpha_0.0_long.txt};

		\addplot+[mark=none,
		draw=none,
       domain=-0.0:2.50,
       samples=100,
       pattern= dots,
       area legend,
       pattern color=darksilver]{-10} \closedcycle;
	
	\nextgroupplot[
		standard,
		title={Reg., $t=0.05$},
 		xlabel={$\diff_t$},
		height=0.22\textwidth,
		width=0.22\textwidth,
        ymin=-0.425,
        ymax=0.425,
        xmin=0.2,
        xmax=1.05,
		ylabel shift=-125pt,
		xtick={0.575},
        xticklabels={\phantom{{$x_1$}, {$x_2$}}},
		yticklabel pos=right,
		xmajorticks=false,
		ymajorticks=false,
		zmajorticks=false,
		ytick={0},
		yticklabels={$0\:$},
		enlarge x limits=0.,
		enlarge y limits=0.,
		view={0}{90},
		colormap name=joshuas,
		]


		\def\t{0.05}
		\def\mval{-0.0089412}
		\def\alf{0.01}
		\addplot3[
    		contour lua,
			contour/levels={0.001, 0.004,0.009,0.016, 0.025, 0.036, 0.049, 0.064, 0.081, 0.1, 0.121, 0.144, 0.169, 0.196, 0.225, 0.256, 0.289, 0.324, 0.361, 0.4, 0.441, 0.484, 0.529, 0.576, 0.625, 0.676, 0.729, 0.784,},
			contour/labels = false,
			samples=75,
			handler/.style=smooth,
			domain = 0.0:1.1,
			thick,
			y domain = 0.0000001:0.5,
		]{(\t ^ 2) * (1.5^2)/ 2 + \mval - \t  * (3/2 * (x - 0.25) - 1 * (y - 0.25)) + 1 / 2 * ((x - 0.25)^2 + (y - 0.25)^2) - \alf * (ln(y / 0.25) - y / 0.25)};

		\draw[darksilver] (axis cs: 0.0,0.0) -- (axis cs: 1.5,0.0);

       	\addplot[joshua,very thick, dashed, lightsilver] table[x index={1},y index={2}, col sep=comma] {csv/regularization/phase_space_alpha_0.01.txt};

       	\addplot[joshua,very thick] table[x index={1},y index={2}, col sep=comma] {csv/regularization/phase_space_alpha_0.01_short.txt};

		\addplot+[mark=none,
		draw=none,
       domain=-0.0:2.50,
       samples=100,
       pattern= dots,
       area legend,
       pattern color=darksilver]{-10} \closedcycle;

	\nextgroupplot[
		standard,
		title={Reg., $t=0.45$},
 		xlabel={$\diff_t$},
		height=0.225\textwidth,
		width=0.225\textwidth,
        ymin=-0.425,
        ymax=0.425,
        xmin=0.2,
        xmax=1.05,
		xtick={0.575},
        xticklabels={\phantom{{$x_1$}, {$x_2$}}},
		yticklabel pos=right,
		xmajorticks=false,
		zmajorticks=false,
		ytick={0},
		yticklabels={$0\:$},
		enlarge x limits=0.,
		enlarge y limits=0.,
		view={0}{90},
		colormap name=joshuas,
		]


		\def\t{0.45}
		\def\mval{0.052583}
		\def\alf{0.01}
		\addplot3[
    		contour lua,
			contour/levels={0.001, 0.004,0.009,0.016, 0.025, 0.036, 0.049, 0.064, 0.081, 0.1, 0.121, 0.144, 0.169, 0.196, 0.225, 0.256, 0.289, 0.324, 0.361, 0.4, 0.441, 0.484, 0.529, 0.576, 0.625, 0.676, 0.729, 0.784,},
			contour/labels = false,
			samples=75,
			handler/.style=smooth,
			domain = 0.0:1.1,
			thick,
			y domain = 0.0000001:0.5,
		]{(\t ^ 2) * (1.5^2)/ 2 + \mval - \t  * (3/2 * (x - 0.25) - 1 * (y - 0.25)) + 1 / 2 * ((x - 0.25)^2 + (y - 0.25)^2) - \alf * (ln(y / 0.25) - y / 0.25)};

		\draw[darksilver] (axis cs: 0.0,0.0) -- (axis cs: 1.5,0.0);

       	\addplot[joshua,very thick, dashed, lightsilver] table[x index={1},y index={2}, col sep=comma] {csv/regularization/phase_space_alpha_0.01.txt};
       	\addplot[joshua,very thick] table[x index={1},y index={2}, col sep=comma] {csv/regularization/phase_space_alpha_0.01_long.txt};

		\addplot+[mark=none,
		draw=none,
       domain=-0.0:2.50,
       samples=100,
       pattern= dots,
       area legend,
       pattern color=darksilver]{-10} \closedcycle;
	\end{groupplot}
\end{tikzpicture}
     \vspace{-6mm}
     \caption{\textbf{Perspective IV: Activation of inequality constraints.} Pressureless Euler equations amount to solving the constrained optimization problem \cref{eqn:opt_problem} (first panel). Shocks form when the constraint becomes active (second panel). 
     Adding a logarithmic barrier function ensures strict feasibility while closely tracking the nominal solution path (third, fourth panel).}
     \label{fig:potentials}
     \vspace{-2mm}
\end{figure}

\subsection{Perspective IV: Optimization}
\label{sec:shock_optimization}  
We conclude by revisiting \cref{eqn:burger_velocity} from an optimization perspective.
In the unidimensional pressureless case and for sufficiently regular $u_0$, the flow map is the minimizer of (\cite{natile2009wasserstein})
\begin{equation}
     \label{eqn:opt_problem}
    \min \limits_{\partial_x \diff \geq 0}  - t \int \limits_{\R} u_0(x) \left(\diff(x) - x\right) \D x + \frac{1}{2}\int \limits_{\R} \left\|\diff(x) - x \right\|^2\D x.
\end{equation}
$\diff_t$ minimizes a linear objective that rewards following the initial velocity, subject to a quadratic regularization. 
As $t$ increases, the relative weight of the regularizer decreases, as in homotopy interior point methods \cite{nemirovski2008interior}.
$\diff_t$ is strictly increasing at first and thus in the interior of the feasible set.
In this regime, it has the closed-form solution $\diff_t(x) = x + t u_0(x)$.
As $t$ increases, $\diff_t$ approaches the boundary of the feasible set, the conic constraint becomes active, and a shock forms (c.f. \cref{fig:potentials}).
As $t$ increases further, $\diff_t$ slides along the boundary of the feasible set.
The first-order optimality conditions of inequality constrained optimization problems require that descent directions point out of the feasible set. 
As $t$ increases, $\diff_t$ pushes against the boundary of the feasible set --- another perspective on the entropy conditions.
\newpage

\section{Information geometric regularization}
\label{sec:infgeoreg}
\subsubsection*{Interior point methods for Euler}
The optimization problem \cref{eqn:opt_problem} resembles homotopy interior point methods of the form
\begin{equation}
     \label{eqn:homotopy}
     \min \limits_{\partial_x \diff \geq 0 } 
          - \int \limits_{\R} u_0(x) \left(\diff(x) - x\right) \D x
          +\frac{1}{t} \underbrace{\left(\bpot(\diff) -\bpot(\cnst{Id}) - \left\langle \nabla\bpot(\cnst{Id}), \left(\diff - \cnst{Id}\right) \right\rangle\right)} \limits_{\Div{\bpot}{\diff}{\mathrm{Id}} \coloneqq },
\end{equation}
with the choice $\bpot_0(\diff) = \int_{\R} \left|\diff(x) - x\right|^2 / 2 \D x$ \cite{nemirovski2008interior}.
These methods solve constrained optimization problems by regularizing them with \emph{Bregman divergences} $\Div{\bpot}{\cdot_1}{\cdot_2}$ generated by \emph{barrier functions} $\bpot$.
Repeatedly solving the problem for decreasing regularization weights $t^{-1}$ yields the \emph{central path} that converges to the solution of the unregularized problem \cite{nemirovski2008interior}.
Choosing a barrier function $\bpot$ that is singular at the boundary of the feasible set enforces strict feasibility of the minimizer for every finite $t$. 
Shock formation in the unidimensional pressureless Euler equation arises from $\bpot_0$ lacking this property. 
This motivates regularizing it by adding the popular log-barrier function \cite{guler1996barrier} to $\bpot$, resulting in the augmented barrier with weight $\alpha \geq 0$, 
\begin{equation}
     \label{eqn:barrier_unidimensional}
     \bpot_{\alpha}(\diff) \coloneqq \int \limits_{\R} \left|\diff(x) - x\right|^2 / 2 \D x + \alpha \int \limits_{\R} - \log(\partial_x \diff(x)) \D x.   
\end{equation}
\cref{fig:potentials} shows that the resulting barrier becomes singular as $\partial_x \diff$ approaches zero and thus produces strictly feasible central paths.
\cref{fig:characteristic_regularization_effort} illustrates how $\alpha$ modulates the regularization and thus the rate at which the central path approaches the boundary of the feasible set.
Translating the central paths into particle trajectories shows that in the regularized system, particles never collide but instead approach each other asymptotically.
Thus, they avoid singularities while maintaining nominal long-time behavior (i.e. in shocks, regularized paths converge to nominal ones as $t \rightarrow \infty$).

\begin{figure}
     \begin{tikzpicture}

\begin{groupplot}[
	compat=1.3,
	group style={group size=4 by 1,
	horizontal sep=0.55cm,
    vertical sep=0.5cm,},
	]
	\nextgroupplot[
		standard,
 		xlabel={$\diff_t(x)$},
		ylabel style={yshift=-6.0pt},	
		ylabel={$t$},
		height=0.200\textwidth,
		width=0.200\textwidth,
        ymin=-0.05,
        ymax=1.0,
        xmin=-0.05,
        xmax=1.0,
        xticklabels={{$x_1$}, {$x_2$}},
        yticklabels={{$0$}},
		xtick={0,0.5},
		ytick={0},
		enlarge x limits=0.,
		enlarge y limits=0.,
		legend style={
			at={(2.20,1.00)},inner sep=3pt,anchor=south,legend columns=5,legend cell align={left}, draw=none,fill=none},
		]	

		\addlegendimage{steelblue,very thick}
		\addlegendimage{joshua,very thick}
		\addlegendimage{orange,very thick}
		\addlegendimage{rust,very thick}
		\addlegendimage{area legend, pattern=dots, pattern color=darksilver}

		\legend{Nominal ($\alpha = 0.0$), $\alpha = 10^{-5}$, $\alpha = 10^{-4}$, $\alpha = 10^{-3}$, infeasible set}

       	\addplot[rust,very thick] table[x index={1},y index={0}, col sep=comma] {csv/regularization/characteristics_alpha_0.001.txt};

        \addplot[rust,very thick] table[x index={2},y index={0}, col sep=comma] {csv/regularization/characteristics_alpha_0.001.txt};
       	\addplot[orange,very thick] table[x index={1},y index={0}, col sep=comma] {csv/regularization/characteristics_alpha_0.0001.txt};

        \addplot[orange,very thick] table[x index={2},y index={0}, col sep=comma] {csv/regularization/characteristics_alpha_0.0001.txt};

       	\addplot[joshua,very thick] table[x index={1},y index={0}, col sep=comma] {csv/regularization/characteristics_alpha_1.0e-5.txt};

        \addplot[joshua,very thick] table[x index={2},y index={0}, col sep=comma] {csv/regularization/characteristics_alpha_1.0e-5.txt};

       	\addplot[steelblue,very thick] table[x index={1},y index={0}, col sep=comma] {csv/regularization/characteristics_alpha_0.0.txt};
        \addplot[steelblue,very thick] table[x index={2},y index={0}, col sep=comma] {csv/regularization/characteristics_alpha_0.0.txt};
		
		\draw[thick, darksilver] (axis cs: 0.575,0.55) rectangle (axis cs: 0.675,0.45);

		\draw[darksilver] (axis cs: 0.0,0.0) -- (axis cs: 1.5,0.0);

		\coordinate (tr1) at (axis cs:0.675,0.55);
		\coordinate (br1) at (axis cs:0.675,00.45);

	\nextgroupplot[
		standard,
 		xlabel={$\diff_t(x)$},
		height=0.200\textwidth,
		width=0.200\textwidth,
        ymin=0.45,
        ymax=0.55,
        xmin=0.575,
        xmax=0.675,
        xticklabels={{$\phantom{x_1}$}},
        yticklabels={},
		xtick={0.575},
		ymajorticks=false,
		enlarge x limits=0.,
		enlarge y limits=0.,
		legend style={
			at={(2.27,1.02)},inner sep=3pt,anchor=south,legend columns=4,legend cell align={left}, draw=none,fill=none},
		]	

       	\addplot[rust,very thick] table[x index={1},y index={0}, col sep=comma] {csv/regularization/characteristics_alpha_0.001.txt};

        \addplot[rust,very thick] table[x index={2},y index={0}, col sep=comma] {csv/regularization/characteristics_alpha_0.001.txt};

       	\addplot[orange,very thick] table[x index={1},y index={0}, col sep=comma] {csv/regularization/characteristics_alpha_0.0001.txt};

        \addplot[orange,very thick] table[x index={2},y index={0}, col sep=comma] {csv/regularization/characteristics_alpha_0.0001.txt};

       	\addplot[joshua,very thick] table[x index={1},y index={0}, col sep=comma] {csv/regularization/characteristics_alpha_1.0e-5.txt};

        \addplot[joshua,very thick] table[x index={2},y index={0}, col sep=comma] {csv/regularization/characteristics_alpha_1.0e-5.txt};

       	\addplot[steelblue,very thick] table[x index={1},y index={0}, col sep=comma] {csv/regularization/characteristics_alpha_0.0.txt};

        \addplot[steelblue,very thick] table[x index={2},y index={0}, col sep=comma] {csv/regularization/characteristics_alpha_0.0.txt};

		\coordinate (tl1) at (axis cs:\pgfkeysvalueof{/pgfplots/xmin},\pgfkeysvalueof{/pgfplots/ymax});
		\coordinate (bl1) at (axis cs:\pgfkeysvalueof{/pgfplots/xmin},\pgfkeysvalueof{/pgfplots/ymin});

	\nextgroupplot[
		standard,
 		xlabel={$\diff_t$},
		ylabel={$\partial_x \diff_t$},
		height=0.200\textwidth,
		width=0.200\textwidth,
        ymin=-0.425,
        ymax=0.425,
        xmin=0.2,
        xmax=1.05,
        xticklabels={{$\phantom{x_1}$}},
		xtick={1.05},
        yticklabels={},
		ymajorticks=false,
		enlarge x limits=0.,
		enlarge y limits=0.,
		]	

		\addplot+[mark=none,
        domain=0:1.25,
        samples=100,
        pattern= dots,
        area legend,
        pattern color=darksilver]{-20} \closedcycle;
	      
		\draw[darksilver] (axis cs: 0.0,0.0) -- (axis cs: 1.5,0.0);

       	\addplot[rust,very thick] table[x index={1},y index={2}, col sep=comma] {csv/regularization/phase_space_alpha_0.001.txt};

       	\addplot[orange,very thick] table[x index={1},y index={2}, col sep=comma] {csv/regularization/phase_space_alpha_0.0001.txt};

       	\addplot[joshua,very thick] table[x index={1},y index={2}, col sep=comma] {csv/regularization/phase_space_alpha_1.0e-5.txt};

       	\addplot[steelblue,very thick] table[x index={1},y index={2}, col sep=comma] {csv/regularization/phase_space_alpha_0.0.txt};

		\draw[thick, darksilver] (axis cs: 0.575,0.05) rectangle (axis cs: 0.675,-0.05);

		\coordinate (tr2) at (axis cs:0.675,0.05);
		\coordinate (br2) at (axis cs:0.675,-0.05);

	\nextgroupplot[
		standard,
 		xlabel={$\diff_t$},
		height=0.200\textwidth,
		width=0.200\textwidth,
        ymin=-0.05,
        ymax=0.05,
        xmin=0.575,
        xmax=0.675,
        xticklabels={{$\phantom{x_1}$}},
		xtick={0.675},
		ytick={0},
		yticklabels={$0\:$},
		enlarge x limits=0.,
		enlarge y limits=0.,
		scaled y ticks=false,
		]	

	 	\addplot+[mark=none,
        domain=0.25:0.75,
        samples=100,
        pattern=dots,
        pattern color=darksilver]{-5} \closedcycle;

		\draw[darksilver] (axis cs: 0.0,0.0) -- (axis cs: 1.5,0.0);

       	\addplot[rust,very thick] table[x index={1},y index={2}, col sep=comma] {csv/regularization/phase_space_alpha_0.001.txt};
	
       	\addplot[orange,very thick] table[x index={1},y index={2}, col sep=comma] {csv/regularization/phase_space_alpha_0.0001.txt};

       	\addplot[joshua,very thick] table[x index={1},y index={2}, col sep=comma] {csv/regularization/phase_space_alpha_1.0e-5.txt};

       	\addplot[steelblue,very thick] table[x index={1},y index={2}, col sep=comma] {csv/regularization/phase_space_alpha_0.0.txt};

		\coordinate (tl2) at (axis cs:\pgfkeysvalueof{/pgfplots/xmin},\pgfkeysvalueof{/pgfplots/ymax});
		\coordinate (bl2) at (axis cs:\pgfkeysvalueof{/pgfplots/xmin},\pgfkeysvalueof{/pgfplots/ymin});

	\end{groupplot}

	\draw[thick, darksilver] (tr2) -- (tl2);
	\draw[thick, darksilver] (br2) -- (bl2);

	\draw[thick, darksilver] (tr1) -- (tl1);
	\draw[thick, darksilver] (br1) -- (bl1);
\end{tikzpicture}  
     \vspace{-8mm}
     \caption{\textbf{Geometric regularization.} Regularizing the geometry of flow maps with $\alpha(-\log(\partial_x \diff))$ avoids singularities and preserves the long-time behavior.}
     \label{fig:characteristic_regularization_effort}
\end{figure}

\subsubsection*{Information geometric regularization}
The variational form \cref{eqn:opt_problem} is limited to the unidimensional pressureless case. 
In general, Euler-type equations are merely dynamical systems on diffeomorphism manifolds. 
To extend \cref{fig:characteristic_regularization_effort} to this case, we exploit the information-geometric interpretation of interior point methods \cite{amari2000methods,amari2016information,ay2017information} described in \cite[Chapter 13.5]{amari2016information}.
Interpreting $\bpot$ as the entropy of a family of probability distributions, the central path of homotopy interior point method with barrier $\bpot$, such as \cref{eqn:homotopy}, is a \emph{dual geodesic}.
Dual geodesics are ordinary (Euclidean) straight lines in the space of gradients $\nabla \bpot$ of the barrier, with general form $t \mapsto \nabla\bpot^{-1}\left( \theta_0 + t v \right)$. 
The gradients of the barrier $\bpot_{\alpha}$ have singularities at the boundary of the feasible set.
Thus, their dual geodesics never reach the boundary (c.f. \cref{fig:characteristic_regularization_effort}).
Information geometric regularization replaces Euclidean geodesics \cref{eqn:characteristic_solution} with dual geodesics generated by $\bpot_{\alpha}$.

\subsubsection*{The multidimensional case}
Positive $\partial_x \diff$ for equations on $\R$ extends to positive Jacobian determinant $\det (\cnst{D} \diff)$ for equations on $\R^d$, suggesting the barrier 
\begin{equation}
     \label{eqn:mod_logdet_barrier}
     \hat{\psi}_{\alpha}(\diff) \coloneqq \int \limits_{\R^d} \left\|\diff(x) - x\right\|^2 / 2 \D x + \alpha \int \limits_{\R^d} (- \log \det(\cnst{D} \diff)) \D x,
\end{equation}
mirroring the transition from linear to positive-semidefinite programming.
Since $-\log \det$ is convex only on symmetric-positive matrices, $\hat{\bpot}$ is non-convex and does not generate a dual geodesic in the sense of \cite{amari2010information}. 
The deeper reason for this problem is that in dimension $d > 1$, connected components of diffeomorphism manifolds are not convex. 
This problem has a long history in nonlinear elasticity, where it prevented the use of convexity in proving the existence of minimizers for realistic energy functionals. 
Ball \cite{ball1976convexity} overcame this problem by introducing so-called polyconvex energy functionals that are convex functions of minors of the gradient.
Motivated by Ball's polyconvexity, we view multivariate diffeomorphism manifolds as a submanifold of the space $\extended$ of pairs $(\diff, \deriv{\diff})$ of flow maps and (formally independent) Jacobian determinants, via the embedding $\diff \mapsto \emb(\diff) \coloneqq (\diff, \det(\cnst{D}\diff))$.
\begin{equation}
     \psi_\alpha(\diff, \deriv{\diff})  \coloneqq \int \limits_{\R^d} \left\|\diff(x) - x\right\|^2 / 2 \D x + \alpha \int \limits_{\R^d} - \log \left(\deriv{\diff}\right) \D x
\end{equation}
is a strongly convex barrier on $\extended$.
The inner product defined by the Hessian of $\bpot$ allows deriving an induced dual geodesic on the diffeomorphism submanifold from that on the ambient $(\diff, \deriv{\diff})$-space, enabling us to treat non-convex problems.

\subsubsection*{Pressure and external forces}
The pressureless Euler equations are geodesic equations on the diffeomorphism manifold. 
More complicated equations like the compressible Euler and Navier-Stokes equations arise by adding forces to the geodesic equations \cite{khesin2021geometric} due to, for instance, pressure, viscosity, and external forces.
Adding these forces to dual geodesic equations yields regularizations of the Euler and Navier Stokes equations and numerous other systems including those listed in \cite[Table 1]{khesin2021geometric}. 
For the sake of concreteness and simpler exposition, this work focuses on the case of the barotropic Euler equation.
For nonzero pressure, shock waves of entropy solutions have finite density and are thus not on the boundary of the diffeomorphism manifold.
But jump discontinuities in the velocity cause the flow to leave the manifold instantaneously. 
Thus, the geodesically complete geometry defined by $\psi$ also prevents the formation of jump discontinuities.

\subsubsection*{Eulerian description}
The final regularization \cref{eqn:reg_euler_intro} arises from converting the Lagrangian representation in terms of the deformation map $\diff$ to an Eulerian description in terms of $\rho$ and $\rho \vct{u}$ by $\vct{u}\left(\diff(x)\right) \coloneqq \dot{\diff}(x)$ and $\rho(\diff(x)) \coloneqq \left(\det\left(\cnst{D}\diff\right)\right)^{-1}$.
Crucially, this equation can be cast in conservation form, enabling the use of discretely conservative numerics that ensure the correct shock propagation speed \cite[Section 12.3]{leveque1992numerical}. 

\subsubsection*{\emph{Information} geometry?\nopunct}
Information geometry views $\bpot$ as the entropy of a family of probability distributions.
Instead, our derivation of $\bpot$ viewed it as the barrier function of a homotopy interior point method.
To reconcile these perspectives, we observe that $\int_{\R^d} - \log \det(\cnst{D}\diff) \D x = \int_{\R^d} \log(\rho) \rho \D x$ is the negative Shannon entropy of the mass density $\rho$.
The Euclidean geometry of $\bpot_0$ on the other hand corresponds to the Wasserstein distance of the pushforward \cite{benamou2000computational}.
Wasserstein geometry, capturing the motion of individual fluid particles, is insensitive to singularity formation.
Information geometry, capturing statistical inference of particle distributions, is highly sensitive to it. 
Thus, from the perspective of \cref{sec:shock_transport}, our regularization prevents singularity formation by blending Wasserstein with information geometry.  
This amounts to viewing $\rho$ as a hybrid between a physical object (Wasserstein geometry) and statistical estimator (information geometry). 
As described in \cref{sec:igm}, this motivates the broader research program of information geometric mechanics.

\section{Derivation in the unidimensional case, using flat geometry}
\label{sec:derivation_flat_univariate}
\subsection*{Notation and Setting} 
We begin by deriving unidimensional information geometric regularization for general barrier functions $\psi$. 
A formal derivation is sufficient to obtain the regularized PDE. 
Thus, we will not discuss the technical complications that arise in the treatment of infinite-dimensional manifolds \cite{schmeding_2022}. 
We consider problems on $\R$ decaying at infinity and postpone the treatment of bounded domains and boundary values to future work.
To this end, we use the Schwartz space $S$ of smooth functions with all derivatives decaying superpolynomially at infinity. 
Using the ambient space $\ambient \coloneqq \Id + S\left(\R;\R\right)$ with $\Id$ the identity function on $\R$, $\manif \subset \ambient$ denotes the diffeomorphism manifold.
Capital Greek letters denote points on $\manif$ and capital Roman letters denote tangent vectors. 

\subsection{Derivation for abstract barrier \texorpdfstring{$\psi$}{}}

\subsubsection{Exponential maps} 
Differential geometry distinguishes \emph{points} $\diff \in \manif$ on a manifold from \emph{directions} $U \in \Tangent{\diff}{\manif}$ in the tangent space of $\manif$ at $\diff$. 
These two notions are connected by the \emph{exponential map} $\Exp{\diff}{U} \in \manif$ that returns a new point in $\manif$ obtained by traveling, for a unit time interval, with initial velocity $U$, starting in $\diff$. 
For instance, we can use the Euclidean structure of $\ambient$ to define the Exponential map $\Exp{\diff}{U} = \diff + U$.
But any given manifold admits a multitude of exponential maps corresponding to different ways of continuing an initial velocity $U$ to a path $t \mapsto \Exp{\diff}{tU}$. 
Paths of this form are called geodesics.  
Riemannian manifolds possess a so-called Levi-Civita exponential map generating locally length-minimizing geodesics.
But geodesics of general exponential maps need not be length-minimizing.

\subsubsection{The dual exponential map} A strictly convex barrier $\psi$ defines a \emph{dual exponential map} on $\manif$ given by 
\begin{equation}
     \DExp{\psi}{ \diff}{U} \coloneqq \nabla \psi^{-1}\Big(\nabla \psi(\diff) + \left[\cnst{D}^2\psi\left(\diff\right)\right]U \Big).
\end{equation}
As shown in \cref{fig:dual_exp}, $\DExp{\psi}{\diff}{U}$ is obtained from the Euclidean exponential map in the coordinate system given by $\nabla \psi$. 
It plays an important role in mirror descent \cite{nemirovskij1983problem,raskutti2015information,schafer2020competitive} and information geometry \cite{amari2000methods}. 
The dual exponential map on the probability simplex induced by the negative Shannon entropy $\psi(p) = \sum_{i = 1}^n \log(p_i) p_i$ arises as the unique exponential map that is invariant under sufficient statistics, geodesically complete, and flat \cite{cencov2000statistical}.
Many other barriers and their geometries arise from parametric families of probability distributions that are (at least formally) submanifolds of infinite-dimensional probability simplices \cite{amari2000methods,fujiwara2022hommage}. 
For instance, $\psi(\mtx{\Theta}) \coloneqq \log \det \mtx{\Theta}$ is the entropy of the family $\{\mathcal{N}(0, \mtx{\Theta}): \mtx{\Theta} \text{ is s.p.d.}\}$ of centered Gaussian distributions. 

\begin{remark}
     In this work, we follow the common practice of considering the gradient $\nabla$ rather than the differential $\cnst{D}$ for defining dual geodesics. 
     This is because the functional-valued differentials appearing in our setting are most conveniently written as inner products with suitable functions. 
     But we emphasize that the dual exponential map can equally well be formulated in the dual space by replacing $\nabla \bpot$ with $\cnst{D} \bpot$. 
\end{remark}

\begin{figure}
    \centering 
    \begin{tikzpicture}

\node[] at (3.60,1.50) (gradmap) {\Large $\xrightarrow[]{\phantom{(}\nabla \bpot\phantom{)^{-1}}}$};
\node[] at (8.55,1.50) (gradmap) {\Large $\xrightarrow[]{(\nabla \bpot)^{-1}}$};
\node[] at (6.25,3.10) (expformula) {\centering ${\color{darksilver} \DExp{{\color{black}\psi}}{{\color{orange} \diff}}{{\color{steelblue} U}}} \coloneqq {\color{darksilver} \nabla \psi^{-1}\Big(}{\color{rust}\nabla \psi(\diff)} {\color{seagreen} +} {\color{darksky} \left[\cnst{D}^2\psi\left(\diff\right)\right]U} {\color{darksilver}\Big)}$};
\begin{groupplot}[
	compat=1.8,
	group style={group size=4 by 1,
	horizontal sep=2.25cm,
    vertical sep=0.5cm,},
	]

	\nextgroupplot[
		standard,
 		xlabel={$\diff$},
		ylabel={$\partial_x \diff$},
		height=0.2025\textwidth,
		width=0.2025\textwidth,
        ymin=-0.10,
        ymax=0.4,
        xmin=0.20,
        xmax=0.6,
		xmajorticks=false,
		ymajorticks=false,
		ytick={},
		yticklabels={$0\:$},
		enlarge x limits=0.,
		enlarge y limits=0.,
		scaled y ticks=false,
		]	

	 	\addplot+[mark=none,
        domain=0.15:1.25,
        samples=100,
        pattern=dots,
        pattern color=seagreen]{-5} \closedcycle;

		\draw[darksilver] (axis cs: 0.0,0.0) -- (axis cs: 1.5,0.0);
		
		\node[draw,circle, orange, very thick] at (axis cs: 0.25,0.25) (initpos) {};
		\draw[steelblue, ultra thick, -latex] (axis cs: 0.25,0.25) -- (axis cs: 0.325,0.2);

	\nextgroupplot[
		standard,
 		xlabel={$\partial_{\diff} \psi(\diff)$},
		ylabel={$\partial_{(\partial_x\diff)} \psi(\diff)$},
		height=0.2025\textwidth,
		width=0.2025\textwidth,
        ymin=-0.70,
        ymax=-0.1,
        xmin=0.10,
        xmax=0.70,
		xmajorticks=false,
		ymajorticks=false,
		ytick={},
		yticklabels={$0\:$},
		enlarge x limits=0.,
		enlarge y limits=0.,
		scaled y ticks=false,
		]

		\node[draw,circle,rust,very thick] at (axis cs: 0.25,-0.15) (initpos) {};

		\draw[darksky, ultra thick, -latex] (axis cs: 0.25,-0.15) -- (axis cs: 0.325,-0.28);

		\node[draw,circle,seagreen,very thick] at (axis cs: 0.51175,-0.6037) (finalpos) {};
		\addplot[seagreen,very thick,dashed] table[x index={1},y index={2}, col sep=comma] {csv/regularization/dual_phase_space_alpha_0.1.txt};

	\nextgroupplot[
		standard,
 		xlabel={$\diff$},
		ylabel={$\partial_x \diff$},
		height=0.2025\textwidth,
		width=0.2025\textwidth,
        ymin=-0.10,
        ymax=0.4,
        xmin=0.20,
        xmax=0.6,
		xmajorticks=false,
		ymajorticks=false,
		ytick={},
		yticklabels={$0\:$},
		enlarge x limits=0.,
		enlarge y limits=0.,
		scaled y ticks=false,
		]	

	 	\addplot+[mark=none,
        domain=0.15:1.25,
        samples=100,
        pattern=dots,
        pattern color=seagreen]{-5} \closedcycle;

		\draw[darksilver] (axis cs: 0.0,0.0) -- (axis cs: 1.5,0.0);
		
		\node[draw,circle, orange, very thick] at (axis cs: 0.25,0.25) (initpos) {};
		\draw[steelblue, ultra thick, -latex] (axis cs: 0.25,0.25) -- (axis cs: 0.325,0.2);


		\node[draw,circle, darksilver, very thick] at (axis cs: 0.51175,0.135315) (dfinalpos) {};

       	\addplot[darksilver,very thick,dashed] table[x index={1},y index={2}, col sep=comma] {csv/regularization/phase_space_alpha_0.1.txt};

	\end{groupplot}

\end{tikzpicture}
    \vspace{-4mm}
    \caption{The barrier function $\psi$ defines a dual exponential map on $\manif$ through the Euclidean exponential map in the coordinate system given by $\nabla \psi$.}
    \label{fig:dual_exp}
\end{figure}

\subsubsection{The dual equations of motion}
Solution paths of interior point methods with barrier $\psi$ are geodesics of $\mathrm{Exp}^{\psi}$. 
We extend interior point methods to dynamical systems by replacing Euclidean geodesics 
\begin{equation}
     \diff_t \coloneqq \Exp{\diff_0}{t\dot{\diff}_0} \quad \text{with dual geodesics} \quad \diff_t \coloneqq \DExp{\psi}{\diff_0}{t\dot{\diff}_0}.
\end{equation}
We still need to account for forces due to, for instance, pressure or gravity.
To this end, we take the second time derivative to deduce the equations of dual inertial motion 
\begin{equation}
     \label{eqn:dual_inertial}
     \ddot{\diff}_t 
+\left[\cnst{D}^2\bpot\left(\diff_t\right)\right]^{-1}\left[\cnst{D}^3\bpot\left(\diff_t\right)\right]\left(\dot{\diff}_t,\dot{\diff}_t\right) = 0.
\end{equation}
In the Euclidean case, $\psi$ is quadratic and the curvature term $\cnst{D}^3\psi$ vanishes, yielding Newton's first law.
For forces $K\big(t, \diff_t, \dot{\diff}_t\big)$ the dual version of Newton's first law is
\begin{equation}
     \label{eqn:second_newton}
     \ddot{\diff}_t 
+\left[\cnst{D}^2\bpot\left(\diff_t\right)\right]^{-1}\left[\cnst{D}^3\bpot\left(\diff_t\right)\right]\left(\dot{\diff}_t,\dot{\diff}_t\right) = K\left(t, \diff_t, \dot{\diff}_t\right).
\end{equation}

\subsection{Deriving the regularized PDE in the unidimensional case} 
\subsubsection{The logarithmic barrier function}
We regularize the squared Euclidean norm with the negative logarithm applied to the derivatives of $\diff$, yielding
\begin{equation}
  \label{eqn:logdetbarrier}
  \bpot(\diff) = \bpot_\alpha(\diff) =  \frac{1}{2} \int \limits_{\R} \abs{\diff(x) - x}^2 \cnst{d}x  + \alpha \int \limits_{\R} -\log \partial_x \diff(x)  \cnst{d}x.
\end{equation}

The parameter $\alpha$ determines the strength of the regularization. 
The negative logarithm is a popular barrier function for optimization on the positive cone. 
It enforces positivity of the iterates and its \emph{self-concordance} ensures rapid convergence of damped Newton minimization algorithms \cite{nemirovski2008interior}. 
It is also the negative Shannon entropy of the Eulerian mass density $\rho(\diff(x)) \defeq \left(\det\left(\cnst{D} \diff(x)\right)\right)^{-1}$. 
Thus, it induces the unique Markov-invariant geometry on the manifold of mass distributions \cite{cencov2000statistical,fujiwara2022hommage}. 
We believe that these connections of $\psi$ to mathematical programming and information geometry will aid the analysis of information geometric regularization in future work. \cite{cao2024information} is a first example.

\subsubsection{The modified equation in Lagrangian coordinates} 
We derive the modified equation by specializing \cref{eqn:second_newton} to the barrier function $\psi$ of \cref{eqn:logdetbarrier}.
This amounts to a function-valued version of computing
\begin{equation}
\left(\frac{(x - 1)^2}{2} - \alpha \log(x)\right)^{\prime \prime \prime} = \left((x - 1) - \alpha \frac{1}{x}\right)^{\prime \prime} = \left(1 + \alpha \frac{1}{x^{2}}\right)^{\prime} = \left(- \alpha \frac{2}{x^{3}}\right).
\end{equation}
Taking the variation in direction $U$ at $\diff$ yields
\begin{equation}
     [\cnst{D} \bpot(\diff)](U) =  \int \limits_{\R} (\diff(x) - x) U(x) \D x + \alpha \int \limits_{\R} - \frac{\partial_x U(x)}{\partial_x \diff(x)} \D x
\end{equation} 
Computing the gradient with respect to the $L^2$ inner product yields 
\begin{equation}
     \nabla \bpot(\diff) = \diff - x +  \alpha  \partial_x \left( \frac{1}{\partial_x \diff} \right).
\end{equation}
Likewise, the second variation is given by  
\begin{equation}
     \left[\cnst{D}^2 \bpot(\diff) \right] \left(U, V\right) 
     = \int \limits_{\R} U(x) \cdot V(x) \D x + \alpha \int \limits_{\R} \frac{ \partial_x U(x) \cdot \partial_x V(x)}{\left(\partial_x \diff(x)\right)^2} \D x,
\end{equation}
and by using the $L^2$ inner product we obtain the linear map
\begin{equation}
     \left[\cnst{D}^2 \bpot(\diff) \right]\colon 
     U \mapsto 
     U
    - \alpha \cdot \partial_x \left(\frac{ \partial_x U}{\left(\partial_x \diff\right)^2}\right).
\end{equation}
Finally, we compute the third variation as 
\begin{equation}
     \left[\cnst{D}^3 \bpot(\diff) \right]
     \left(U, V, W\right) 
     =
     - 2 \alpha \int \limits_{\R} \frac{\partial_x U(x) \cdot \partial_x V(x) \cdot \partial_x W(x) }{\left(\partial_x \diff(x)\right)^3} \D  x.
\end{equation}
Again using the $L^2$ inner product, we obtain the quadratic map 
\begin{equation}
     \left[\cnst{D}^3 \bpot(\diff) \right]
     \colon
     \left(U, V\right) 
     \mapsto 
     2\alpha \partial_x \left(\frac{\partial_x U \cdot \partial_x V}{(\partial_x \diff)^3}\right)
\end{equation}
Plugging the above into \cref{eqn:second_newton} results in the modified equation 
\begin{equation}
     \label{eqn:lagrangian_eqn}
     \ddot{\diff}_t 
     +\left( (\cdot) - \alpha \partial_x  \left(\left[\partial_x \diff_t \right]^{-2} [\partial_x (\cdot)]  \right) \right)^{-1}  
     2 \alpha \partial_x \left(\left[\partial_x \diff_t \right]^{-3}  [\partial_x \dot{\diff}_t]^2 \right) = K\left(t, \diff_t, \dot{\diff}_t\right).
\end{equation}
Equivalently we can write
\begin{equation}
     \label{eqn:lagrangian_eqn_v2}
     \left((\cdot) - \alpha \partial_x  \left(\left[\partial_x \diff_t \right]^{-2} [\partial_x (\cdot)] \right) \right) \left(\ddot{\diff}_t - K\left(t, \diff_t, \dot{\diff}_t\right)\right)
     = - 2 \alpha \partial_x \left(\left[\partial_x \diff_t \right]^{-3}  [\partial_x \dot{\diff}_t]^2 \right). 
\end{equation}

\subsubsection{The modified equations in Eulerian coordinates}
We write the Eulerian velocity field $u$ as 
\begin{equation}
     \dot{\diff}_t(x) \eqqcolon u\left(\diff_t(x)\right) \quad \Rightarrow \quad  
     \ddot{\diff}_t(x) = \left[\partial_x u (\diff_t(x))\right] u\left(\diff_t\left(x\right)\right)  + \partial_t u(\diff_t(x))
\end{equation}
and the Eulerian density $\rho$ as $\rho(\diff_t(x)) = \left(\partial_x \diff_t(x)\right)^{-1}$.
For $V \coloneqq \ddot{\diff}_t - K(t, \diff_t, \dot{\diff}_t)$ we integrate both sides of \cref{eqn:lagrangian_eqn_v2} against a vector field $U$, obtaining
\begin{equation}
     \int \limits_{\R } U \cdot \left(V - \alpha \partial_x  \left(\left[\partial_x \diff_t \right]^{-2} [\partial_x V]  \right) \right)  \D x 
     =- \int \limits_{\R } 2 \alpha U \cdot \partial_x \left(\left[\partial_x \diff_t \right]^{-3}  [\partial_x \dot{\diff}_t]^2 \right) \D x.
\end{equation}
Integrating by parts, we obtain 
\begin{equation}
     \int \limits_{\R} U \cdot V + \alpha \left[\partial_x U\right] \left[\partial_x \diff_t \right]^{-2} [\partial_x V]  \D x 
     = \int \limits_{\R} 2 \alpha \left[\partial_x U\right] \left[\partial_x \diff_t \right]^{-3}  [\partial_x \dot{\diff}_t]^2
      \D x.
\end{equation}
We now perform a change of variables 
\begin{align}
     &\int \limits_{\R} \left(U \cdot V
     + \alpha  \left[\partial_x U\right] \left[\partial_x \diff_t \right]^{-2} [\partial_x V]  \right) \circ \diff_{t}^{-1} \cdot \left(\partial_x \diff_{t}^{-1}\right) \D x \\
     =& \int \limits_{\R} 2 \alpha \big(\left[\partial_x U\right] \left[\partial_x \diff_t \right]^{-3}  [\partial_x \dot{\diff}_t]^2  \big) \circ \diff_t^{-1} \cdot \left(\partial_x \diff_{t}^{-1}\right)\D x.
\end{align}
Noting that $[\partial_x \dot{\diff}_t] = [\partial_x u] \circ \diff_{t} [\partial_x \diff_t]$ and writing $U = \tilde{U} \circ \diff$ and $V = \tilde{V} \circ \diff$ yields
\begin{equation}
     \int \limits_{\R} \left(\tilde{U} \cdot \tilde{V} + \alpha  [\partial_x \tilde{U}] [\partial_x \tilde{V}] \right)  \rho \D x 
     = \int \limits_{\R} 2 \alpha [\partial_x \tilde{U}]  [\partial_x u]^2  \rho \D x.
\end{equation}
Integrating by part once more, we obtain as before
\begin{equation}
     \int \limits_{\R} \tilde{U} \cdot \tilde{V} \rho - \tilde{U} \cdot \alpha \partial_x \left( [\partial_x \tilde{V}] \rho \right)  \D x 
     = -\int \limits_{\R} 2 \alpha \tilde{U} \cdot \partial_x \left( [\partial_x u]^2 \rho \right) \D x.
\end{equation}
Using the fact that $U$ and thus $\tilde{U}$ was chosen fully general, we obtain 
\begin{equation}
     \partial_t u + [\partial_x u] u + \left((\cdot) \rho - \alpha \partial_x \left(\rho [\partial_x (\cdot)]\right)\right)^{-1} 2 \alpha \partial_x \left(\rho [\partial_x u]^2 \right) = \rho K(t, \diff_t, \dot{\diff}_t) \circ \diff_{t}^{-1},
\end{equation}
which we combine with the transport equation $\partial_t \rho + \partial_x \left(u \rho\right) = 0$.

\subsubsection{A regularized conservation law}
\label{sec:conservation_1d}
To express the regularized equation in conservation form, we multiply with $\rho$ and use the equation for mass conservation,  
\begin{align}
     \partial_t \left(\rho u\right) + \left[\partial_x \rho \right] u^2 + u \rho (\partial_x u) + \rho \left[\partial_x u\right] u
           + 2 \alpha \rho \left( \rho - \alpha \partial_x \rho \partial_x \right)^{-1} \partial_x \left( \rho \left[\partial_x u \right]^2 \right)
     \! = \! \rho K \circ \diff_{t}^{-1} \nonumber.
\end{align}
We now use the identity $\partial_x \left( \rho u \cdot u \right) = \partial_x (\rho u) u + \rho u \partial_x (u)$
to conclude 
\begin{equation}
     \partial_t \left(\rho u\right) + \partial_x\left(\rho u^2\right)
           + 2 \alpha \rho \left( \rho - \alpha \partial_x \rho \partial_x \right)^{-1} \partial_x \left( \rho \left[\partial_x u\right]^2\right) = \rho K(t, \diff_t, \dot{\diff}_t) \circ \diff_{t}^{-1}.
\end{equation}
We now use the identity $\left(\mathcal{A} \mathcal{B}\right)^{-1} = \mathcal{B}^{-1} \mathcal{A}^{-1}$ for linear operators $\mathcal{A}$ and $\mathcal{B}$ to obtain
\begin{align}
       & 2 \alpha \rho \left( \rho - \alpha \partial_x \rho \partial_x \right)^{-1} \partial_x \left( \rho \left[\partial_x u\right]^2\right)
     =  2 \alpha  \left( \mathrm{Id} - \alpha \partial_x \rho \partial_x \rho^{-1}\right)^{-1} \partial_x \left( \rho \left[\partial_x u\right]^2\right)\\
     = & 2 \alpha  \left( \partial_x^{-1} - \alpha \rho \partial_x  \rho^{-1}\right)^{-1} \left(\rho \left[\partial_x u\right]^2\right)
     = 2 \alpha  \partial_x \left( \left( \mathrm{Id} - \alpha \rho \partial_x \rho^{-1} \partial_x \right)^{-1} \left(\rho \left[\partial_x u \right]^2\right)\right)\\
     = & 2 \alpha  \partial_x \left( \left( \rho^{-1} - \alpha \partial_x \rho^{-1} \partial_x \right)^{-1} \left[\partial_x u\right]^2\right) \eqqcolon \partial_x(\Sigma).
\end{align}
Here, $\rho, \rho^{-1}$ inside $\left(\cdot\right)^{-1}$ are understood as the linear multiplication-by-$\rho, \rho^{-1}$ operators. 
The inverse derivative $\partial_x^{-1}$ is understood as a mapping from $S$ to $S$. 
Setting $ K(t, \diff_t, \dot{\diff}_t) \circ \diff_{t}^{-1} = \partial_{x} P(\rho) - f$, the regularized Euler equation becomes 
\begin{equation}
     \label{eqn:infgeo_reg_flux_1d}
     \begin{cases}          
          \partial_t(\rho u) + \partial_x \left(\rho u^2+ P(\rho) + \Sigma  \right) &= f\\
          \partial_t \rho + \partial_x(\rho u) &= 0\\
          \Sigma \rho^{-1} - \alpha \partial_x \left(\rho^{-1} \partial_x \Sigma \right) &=  2 \alpha [\partial_x u]^2.
     \end{cases}
\end{equation}

\section{Derivation in the multivariate case, using a nonlinear embedding}
\label{sec:derivation_multivariate}
\subsection{The need for non-convex information geometry}
\subsubsection{Shock formation in multidimensional problems} In the univariate setting, the only scenario for shock formation is the vanishing of the derivative $\partial_x \diff_t$ of the deformation map. 
In $d>1$ dimensions, any loss of invertibility of $\cnst{D} \diff_t$ results in shock formation. 
The dimension of the null space of $\cnst{D} \diff_t$ determines the dimension of the resulting shock singularity. 
The most common situation, a $(d-1)$-dimensional shock front, arises from $\cnst{D} \diff_t$ having a one-dimensional null space normal to the shock front. 
Thus, our barrier function $\psi$ needs to prevent any loss of invertibility of $\cnst{D} \diff_t$. 

\subsubsection{The log determinantal barrier function}
In mathematical programming, the natural generalization of the logarithmic barrier function to semidefinite programming is the log determinant barrier function \cite{guler1996barrier}.  
Just like the logarithmic barrier function on the positive cone, it is convex and self-concordant on the cone of symmetric and positive-definite matrices, enabling the efficient implementation of homotopy methods.
This suggests generalizing \cref{sec:derivation_flat_univariate} by using the barrier
\begin{equation}
     \label{eqn:mod_logdet_barrier_derivation}
     \hat{\bpot}(\diff) = \hat{\bpot}_\alpha(\diff) \coloneqq \frac{1}{2} \int \limits_{\R^d} \left\|\diff - \mathrm{Id}\right\|^2 \D x  + \alpha \int \limits_{\R^d} (- \log \det(\cnst{D} \diff)) \D x.
\end{equation}

\subsubsection{Nonconvexity of $\hat{\bpot}$}
The problem with this approach is that the Jacobians of diffeomorphisms are generally not symmetric and thus $\hat{\bpot}$ is not convex on $\manif$ and does not define a dually flat information geometry in the sense of \cite{amari2000methods}. 
Likewise, in mathematical elasticity, no convex energy functional can simultaneously be rotationally invariant and prevent the collapse of deformation maps \cite{coleman1959thermostatics}. 
The seminal work \cite{ball1976convexity} overcomes this problem by introducing polyconvex energy functionals: convex functions of Jacobian determinants. 
The determinant's weak continuity allows applying the direct method of the calculus of variations to polyconvex problems \cite{dacorogna2007direct}. 
\subsubsection{Regularization by embedding} Inspired by polyconvexity, we introduce $\xi\colon \diff \mapsto \left(\diff, \det\left(\cnst{D} \diff\right)\right)$ that embeds the diffeomorphism manifold $\manif$ into the extended ambient space $\extended$ consisting of pairs $\left(\diff, \deriv\diff\right)$ of deformation and Jacobian-determinant degrees of freedom.
We define the barrier function $\bpot$ on $\extended$ as 
\begin{equation}
     \label{eqn:extended_barrier}
     \bpot(\diff, \deriv \diff) = \bpot_{\alpha}(\diff, \deriv \diff) \coloneqq \frac{1}{2} \int \limits_{\R^d} \left\|\diff - \mathrm{Id}\right\|^2 \D x  + \alpha \int \limits_{\R^d} (- \log ( \deriv \diff)) \D x.
\end{equation}
This barrier function is strictly convex and thus defines a dually flat information geometry on $\extended$. 
The Riemannian metric defined by the Hessian of $\psi$ and the embedding $\xi$ induce an exponential map on $\manif$ that avoids the formation of shock singularities. 
We think that the weak continuity of $\xi$ will aid the rigorous analysis of the resulting PDE in future work, extending the unidimensional results by \cite{cao2024information}. 
\begin{figure}
     \centering
     \begin{tikzpicture}

	\begin{scope}
		\draw [-stealth, darksilver, thick] (0.0, 1.25) -> (4, 1.25);
		\draw [-stealth, darksilver, thick] (1.25, 0.0) -> (1.25, 3);
		\draw [steelblue, thick] (1.5, 1.5) circle (1);
		\draw [-stealth, joshua, dashed, thick] (1.5, 2.5) -> (1.5 + 1.0, 2.5);
		\draw [-stealth, darksky, thick, bend right=50] (1.5, 2.5) to ({1.5 + sin(30)}, {1.5 + cos(30)});
		\draw [-stealth, orange, dotted, thick] ({1.5 + sin(30)}, {1.5 + cos(30)}) -> ({1.0 + 1.5 + sin(30)}, {1.5 + cos(30)});
		\draw [-{stealth}{stealth}, seagreen, thick] ({1.0 + 1.5 + sin(30)}, {1.5 + cos(30)}) -> ({1.5 + sin(61)}, {1.5 + cos(61)});
		\draw [-stealth, rust, thick] ({1.5 + sin(30)}, {1.5 + cos(30)}) -> ({1.5 + sin(61)}, {1.5 + cos(61)});
		\node[darksilver] at (0.5, 2.5) {$\R^2$};
		\node[draw, steelblue] at (3.5, 2.0) {$S^1$};
		\node at (2.5, 1.00) {$x$};
		\node at (2.5, 1.00) {$x$};
		\node at (1.00, 2.0) {$y$};
	\end{scope}

	\begin{scope}[xshift=42.5,yshift=35]
		\node[text width=4.25cm] at (5.3, 0.25) {{\baselineskip=1pt $\emb(X)$ to $\emb(X + \epsilon V)$ \newline tangent vector $U$ at $\xi(X)$ \newline parallel transport in {\color{darksilver}$\R^2, \extended$} \newline projection onto {\color{steelblue} $\Tangent{}{S^1}$, $\Tangent{}{\extended}$}} \newline parallel transport in {\color{steelblue} $S^1$, $\extended$}\par};
		\draw [-stealth, darksky, thick, bend left=50] (2.60, 1.0) to ({3.05}, {1.0});
		\draw [-stealth, joshua, dashed, thick] (2.60, 0.65) to ({3.05}, {0.65});
		\draw [-stealth, orange, dotted, thick] (2.60, 0.20) to ({3.05}, {0.20});
		\draw [-{stealth}{stealth}, seagreen, thick] (2.60, -0.20) to ({3.05}, {-0.20});
		\draw [-stealth, rust, thick] (2.60, -0.6) to ({3.05}, {-0.6});
	\end{scope}

	\begin{scope}[xshift=255]
		\draw [-stealth, darksilver, thick] (0.0, 1.25) -> (4, 1.25);
		\draw [-stealth, darksilver, thick] (1.25, 0.0) -> (1.25, 3);
		\draw [steelblue, thick] (1.25, 1.25) -- (3.0, 3.0);
		\draw [-stealth, joshua, dashed, thick] (2.6, 2.6) -> (2.2, 2.2);
		\draw [-stealth, darksky, thick, bend right=50] (2.6, 2.6) to (2.05, 2.05);
		\draw [-stealth, orange, dotted, thick] (2.05, 2.05) -> (2.15, 1.35);
		\draw [-{stealth}{stealth}, seagreen, thick] (2.15, 1.35) -> (1.75, 1.75);
		\draw [-stealth, rust, thick] (2.05, 2.05) -> (1.75, 1.75);
		\node[darksilver] at (0.5, 2.5) {$\extended$};
		\node[draw, steelblue] at (3.5, 2.0) {$\manif$};
		\node at (2.5, 1.00) {$\diff$};
		\node at (0.50, 2.0) {$\deriv{\diff} = \partial_x \diff$};
	\end{scope}

\end{tikzpicture}
     \vspace{-4mm}
     \caption{\textbf{Parallel transport on submanifolds.} Parallel transport on a submanifold is parallel transport in the ambient manifold, followed by projection on the tangent space of the submanifold.
     For $S^1 \hookrightarrow \R^2$, the difference between parallel transport on ambient and submanifold arises from the nonlinearity of the embedding $\emb$. 
     In unidimensional IGR, $\emb$ is affine but the parallel transport on the ambient space does not respect its linear structure, causing parallel transport on $\manif$ and $\extended$ to differ. 
     In multidimensional IGR, $\emb$ is nonlinear and both of the above effects contribute.}
     \label{fig:pointing_out}
\end{figure}
\subsection{Derivation for general embedding and potential} 
\subsubsection*{Notation and setting}
We now derive information geometric regularization in the multidimensional setting. 
To this end, we consider the extended ambient space $\extended \coloneqq \left(\Id + S\left(\R^d; \R^d\right)\right) \times \left(1 + S\left(\R^d\right)\right)$, denoting as $S$ the Schwartz space $S$ of smooth functions with all derivatives decaying superpolynomially at infinity. 
We consider the diffeomorphism manifold $\manif \hookrightarrow \extended$ as a submanifold of $\extended$, through the embedding $\xi\colon \diff \mapsto \left(\diff, \det\left(\cnst{D}\diff\right)\right)$.
Our potential $\psi$ will be defined on the extended ambient space $\extended$.
We denote by capital Greek letters points on the manifolds $\manif$ and $\extended$ and by capital Roman letters the elements of their tangent spaces. 
We denote elements of $\extended$ or its tangent spaces as $\left(\diff, \deriv{\diff}\right)$ or $\left(U, \deriv{U}\right)$ to distinguish the two degrees of freedom.

\subsubsection{Exponential maps, geodesics, and parallel transport}
In \cref{sec:derivation_flat_univariate}, we derived dual equations of inertial motion \cref{eqn:dual_inertial} by taking derivatives of the dual exponential map. 
Equivalently, they can be defined via the dual \emph{affine connection}
\begin{equation}
    \Gamma_{\diff}^{\hat{\bpot}} \left(U, V\right) \coloneqq \left[\cnst{D}^2\hat{\bpot}\left(\diff\right)\right]^{-1}\left[\cnst{D}^3\hat{\bpot}\left(\diff\right)\right]\left(U,V\right).
\end{equation}
The affine connection $\Gamma^{\hat{\bpot}}_{\diff}\left(U, \cdot\right)$ defines a correspondence between the tangent spaces $\Tangent{\diff}{\manif}$ and $\Tangent{\diff + \epsilon U}{\manif}$ in infinitesimally nearby base points $\diff$ and $\diff + \epsilon U$. 
It identifies the vector $V + \epsilon \Gamma^{\hat{\bpot}}_{\diff}\left(U, V\right)$ as the tangent vector in $\Tangent{\diff}{\manif}$ corresponding to the tangent vector $V$ in $\Tangent{\diff + \epsilon U}{\manif}$.
Equivalently, we can think of the map $V \mapsto V - \epsilon \Gamma^{\hat{\bpot}}_{\diff}\left(U, V\right)$ from $\Tangent{\diff}{\manif}$ to $\Tangent{\diff + \epsilon U}{\manif}$ as transporting a tangent vector $V \in \Tangent{\diff}{\manif}$ an infinitesimal distance in the direction $\epsilon U$, resulting in an element of $\Tangent{\diff + \epsilon U}{\manif}$.
This operation is called \emph{parallel transport}.
Solutions of dual equations of motion or geodesic equations
\begin{equation}
     \ddot{\diff}_t 
    + \Gamma_{\diff}^{\hat{\bpot}} \left(\dot{\diff}_t, \dot{\diff}_t\right) = 0 
\end{equation}
are therefore curves with velocities that are parallel transported along themselves.

\subsubsection{Parallel transport on submanifolds} 
In more than one spatial dimension, the barrier function $\hat{\bpot}$ defined in \cref{eqn:mod_logdet_barrier_derivation} is not convex on $\manif$ and thus does not define an affine connection on $\manif$.
Instead, we use the parallel transport on $\extended$ induced by $\bpot$ (which is always convex) to induce a parallel transport on the submanifold $\manif \hookrightarrow \extended$ through the embedding $\xi$.   
As outlined in \cite[Section 5.10]{amari2016information}, we proceed in three steps to transport a tangent vector $V$ from $\Tangent{\diff}{\manif}$ to $\Tangent{\diff + \epsilon U}{\manif}$.     
We use the embedding $\xi$ to identify $V$ with the element $[\cnst{D}\xi(\diff)] V$ of the tangent space $\Tangent{\xi (\diff)}{\extended}$ of the ambient manifold $\extended$.
We then parallel transport $[\cnst{D}\xi(\diff)] V$ in the ambient manifold $\extended$ to the point $\xi(\diff + \epsilon U)$, obtaining $(W, \deriv W) \in \Tangent{\xi\left(\diff + \epsilon U\right)}{\extended}$.
Finally, we project $(W, \deriv W)$ onto $\Tangent{\diff + \epsilon U}{\manif}$ via the Riemannian metric on $\extended$ given by the Hessian of $\psi$, 
%
\begin{equation}
    \left \llangle \left(U, \deriv{U}\right), \left(V, \deriv{V}\right)\right \rrangle_{\left(\diff, \deriv{\diff}\right)} \coloneqq \left[\cnst{D}^2\bpot (\diff, \deriv{\diff})\right]\left(\left(U, \deriv{U}\right), \left(V, \deriv{V}\right)\right).
\end{equation}
This is the Fisher-Rao metric of a family of probability distributions with entropy $\bpot$, the unique Riemannian metric invariant under Markov isomorphisms \cite{cencov2000statistical,fujiwara2022hommage}. 
\subsubsection{Computing the induced parallel transport} 
Transporting $[\cnst{D}\xi(\diff)]V$ to $\xi(\diff + \epsilon U)$ yields, to leading order in $\epsilon$,  
\begin{equation}
    (W, \deriv W) \approx [\cnst{D}\xi(\diff)]V - \epsilon \left[\cnst{D}^2 \bpot\left(\xi\left(\diff\right)\right)\right]^{-1} \left[\cnst{D}^3\bpot\left(\xi\left(\diff\right)\right)\right]\left([\cnst{D}\xi(\diff)]U, [\cnst{D}\xi(\diff)]V\right).
\end{equation}
The projection $X$ of $(W, \deriv W)$ onto $\Tangent{\diff + \epsilon U}{\manif}$ satisfies, for all $Y \in \Tangent{\diff + \epsilon U}{\manif}$, 
\begin{equation}
    \left \llangle [\cnst{D} \xi\left(\diff + \epsilon U\right)] Y, [\cnst{D} \xi\left(\diff + \epsilon U\right)] X\right \rrangle_{\xi(\diff + \epsilon U)} =
    \left \llangle [\cnst{D} \xi\left(\diff + \epsilon U\right)] Y, \left(W, \deriv{W}\right) \right \rrangle_{\xi(\diff + \epsilon U)}.
\end{equation}
Using the $L^2$ inner product on $T_{(\diff, \deriv{\diff})}\extended$, we can view $\left[\cnst{D}^2 \bpot \left(\diff, \deriv{\diff}\right)\right]$ as a linear map from $T_{(\diff, \deriv{\diff})}\extended$ to itself, with
\begin{equation}
    \left \llangle \left(U, \deriv{U}\right), \left(V, \deriv{V}\right)\right \rrangle_{\left(\diff, \deriv{\diff}\right)} = \left\langle\left(U, \deriv{U}\right), \left[\cnst{D}^2\bpot (\diff, \deriv{\diff})\right] \left(V, \deriv{V}\right)\right\rangle_{L^2}.
\end{equation}
Denoting the $L^2$-adjoint of $[\cnst{D} \xi\left(\diff + \epsilon U\right)]$ as $[\cnst{D} \xi\left(\diff + \epsilon U\right)]^\ad$, we obtain 
\begin{align}
    X \approx & \left([\cnst{D} \xi\left(\diff + \epsilon U\right)]^\ad \left[\cnst{D}^2\bpot \left(\xi\left(\diff + \epsilon U\right)\right)\right] [\cnst{D} \xi\left(\diff + \epsilon U\right)]\right)^{-1} \\
    &\quad [\cnst{D} \xi\left(\diff + \epsilon U\right)]^\ad \left[\cnst{D}^2\bpot \left(\emb \left(\diff + \epsilon U\right)\right)\right]\\ 
    &\quad \quad \left([\cnst{D}\xi(\diff)]V - \epsilon \left[\cnst{D}^2\bpot\left(\xi\left(\diff\right)\right)\right]^{-1} \left[\cnst{D}^3\bpot\left(\xi\left(\diff\right)\right)\right]\left([\cnst{D}\xi(\diff)]U, [\cnst{D}\xi(\diff)]V\right)\right).
\end{align}
To order $\epsilon^0$, we have $X = V$. In order to compute the first-order approximation of $X$ in $\epsilon$, we compute the derivative of the right-hand side. 
To simplify the notation, we define $\beta_{\diff} \coloneqq [\cnst{D} \xi\left(\diff\right)]$ and $\gamma_{\diff} \coloneqq \left[\cnst{D}^2\bpot \left(\xi\left(\diff\right)\right)\right]$.
We compute

\begin{align}
    \frac{\D}{\D \epsilon}\bigg|_{\epsilon = 0}&\left(\beta_{\diff + \epsilon U}^\ad \gamma_{\diff + \epsilon U} \beta_{\diff + \epsilon U}\right)^{-1}
    \beta_{\diff + \epsilon U}^\ad \gamma_{\diff + \epsilon U}
    \left(\beta_\diff V - \epsilon \gamma_{\diff}^{-1} \left[\cnst{D}^3\bpot\left(\xi\left(\diff\right)\right)\right]\left(\beta_{\diff} U, \beta_{\diff} V\right)\right) \\
    =&- \left(\beta_{\diff}^\ad \gamma_{\diff} \beta_{\diff}\right)^{-1} 
         \left[\frac{\D}{\D \epsilon} \bigg|_{\epsilon = 0} \left(\beta_{\diff + \epsilon U}^\ad \gamma_{\diff + \epsilon U} \beta_{\diff + \epsilon U} \right) \right]
         \left(\beta_{\diff}^\ad \gamma_{\diff} \beta_{\diff}\right)^{-1} 
         \beta_{\diff}^\ad \gamma_{\diff}\beta_\diff V\\
    &+ \left(\beta_{\diff}^\ad \gamma_{\diff} \beta_{\diff}\right)^{-1}
    \left[\frac{\D}{\D \epsilon}\bigg|_{\epsilon = 0} \beta_{\diff + \epsilon U}^\ad \gamma_{\diff + \epsilon U} \right]
    \beta_\diff V\\
    &- \left(\beta_{\diff}^\ad \gamma_{\diff}  \beta_{\diff}\right)^{-1}\beta_{\diff}^\ad \gamma_{\diff} 
    \gamma_{\diff}^{-1}
    \left[\cnst{D}^3\bpot\left(\xi\left(\diff\right)\right)\right]\left(\beta_{\diff} U, \beta_{\diff} V\right)\\
    =& - \left(\beta_{\diff}^\ad \gamma_{\diff} \beta_{\diff}\right)^{-1}\left(
    \beta_{\diff}^\ad \gamma_{\diff}
    \left[\cnst{D}^2 \xi\left(\diff\right) \right] \left(U, V\right)
    + \beta_{\diff}^\ad \gamma_{\diff} \gamma_{\diff}^{-1}
    \left[\cnst{D}^3\bpot\left(\xi\left(\diff\right)\right)\right]\left(\beta_{\diff} U, \beta_{\diff} V\right)\right)
\end{align}
In summary, we obtain to first order in $\epsilon$ that
\begin{align}
    X \approx V &- \epsilon \left([\cnst{D} \xi\left(\diff\right)]^{\ad}\left[\cnst{D}^2\bpot \left(\xi\left(\diff\right)\right)\right][\cnst{D}\xi\left(\diff\right)]\right)^{-1}
    [\cnst{D} \xi\left(\diff\right)]^\ad\left[\cnst{D}^2\bpot \left(\xi\left(\diff\right)\right)\right]\\
    & \left(\left[\cnst{D}^2 \xi\left(\diff\right) \right]\left(U, V\right) + \left[\cnst{D}^2 \bpot\left(\xi\left(\diff\right)\right)\right]^{-1}\left[\cnst{D}^3\bpot\left(\xi\left(\diff\right)\right)\right]\left(\left[\cnst{D}\xi\left(\diff\right)\right] U, \left[\cnst{D}\xi\left(\diff\right)\right] V\right)\right).
\end{align}
Rewriting this result in terms of an affine connection, we obtain
\begin{align}
    \Gamma_{\diff}^{\bpot, \emb}(U, &V) =  \left([\cnst{D} \xi\left(\diff\right)]^{\ad}\left[\cnst{D}^2\bpot \left(\xi\left(\diff\right)\right)\right][\cnst{D}\xi\left(\diff\right)]\right)^{-1} 
    [\cnst{D} \xi\left(\diff\right)]^{\ad}\left[\cnst{D}^2\bpot \left(\xi\left(\diff\right)\right)\right]\\
    &\left(\left[\cnst{D}^2 \xi\left(\diff\right) \right]\left(U, V\right) + \left[\cnst{D}^2\bpot\left(\xi\left(\diff\right)\right)\right]^{-1}\left[\cnst{D}^3\bpot\left(\xi\left(\diff\right)\right)\right]\left(\left[\cnst{D}\xi\left(\diff\right)\right] U, \left[\cnst{D}\xi\left(\diff\right)\right] V\right)\right).
\end{align}
Note that the two summands in the above formula correspond to the two separate contributions to the affine connection, as illustrated in \cref{fig:pointing_out}. 
The first summand arises from the nonlinearity (with respect to the vector space structure of $\extended$) of the embedding $\xi$. 
This is akin to the curvature of a $(d - 1)$-sphere arising from the nonlinearity of its embedding into $\R^d$. 
The second summand arises from the dual parallel transport on $\extended$ not respecting the vector space structure of $\extended$. 
It induces a nontrivial parallel transport even for $d = 1$, where $\xi$ is linear. 
In this case, it gives rise to the equations derived in \cref{sec:derivation_flat_univariate}.

\subsection{The geodesic equation in Lagrangian form}

We now plug in $\emb$ and $\bpot$ and compute $\Gamma^{\psi, \emb}$, in order to derive the equation of inertial motion.
The quadratic form associated with $\left(\cnst{D}^2\bpot\right)$ is given by 
\begin{equation}
    [\cnst{D}^2 \bpot(\diff, \deriv{\diff})]\left(\left(U, \deriv{U}\right), \left(V, \deriv{V}\right)\right) = \int \limits_{\R^{d}} \inprod{U, V}  \D x +  \alpha \int \limits_{\R^{d}} \frac{\deriv{U}\deriv{V}}{\left(\deriv{\diff}\right)^2} \D x  
\end{equation}
and using the $L^2$ inner product on $T_{(\diff, \deriv{\diff})}\extended$, we obtain the linear map 
\begin{equation}
    [\cnst{D}^2 \bpot(\diff, \deriv{\diff})]\left(U, \deriv{U}\right) = \left(U, \alpha \frac{\deriv{U}}{\left(\deriv \diff\right)^2}\right).
\end{equation}
Similarly, we obtain the cubic form of the third derivative as 
\begin{equation}
    [\cnst{D}^3 \bpot(\diff, \deriv{\diff})]\left(\left(U, \deriv{U}\right), \left(V, \deriv{V}\right), \left(W, \deriv{W}\right)\right) = -  2 \alpha \int \limits_{\R^{d}} \frac{\deriv{U}\deriv{V}\deriv{W}}{\left(\deriv{\diff}\right)^3} \D x  
\end{equation}
and the resulting quadratic map 
\begin{equation}
    [\cnst{D}^3 \bpot(\diff, \deriv{\diff})]\left(\left(U, \deriv{U}\right), \left(V, \deriv{V}\right)\right) = \left(0, - 2 \alpha \frac{\deriv{U}\deriv{V}}{\left(\deriv{\diff}\right)^3}\right).
\end{equation}
Similarly, using Jacobi's formula, we can compute the map 
\begin{equation}
    \left[\cnst{D} \emb\left(\diff\right) \right]U = \left(U, \det\left(\cnst{D} \diff\right) \trace\left([\cnst{D}\diff]^{-1} [\cnst{D} U]\right) \right). 
\end{equation}
To compute its adjoint, we multiply with $(W, \deriv{W})$ and integrate 
\begin{equation}
    \langle (W, \deriv{W}), \left[\cnst{D} \emb\left(\diff\right) \right]U \rangle_{L^2} = \int\limits_{\R^d} \inprod{U, W} \D x + \int\limits_{\R^d} \deriv{W} \det\left(\cnst{D} \diff\right) \trace\left([\cnst{D}\diff]^{-1} [\cnst{D} U]\right) \D x. 
\end{equation}
Letting $\overline{V} \coloneqq [\cnst{D} \emb(\diff)]^{\ad} (W, \deriv{W})$, we have, for $\cdiv$ the column-wise divergence,
\begin{align}
    \int\limits_{\R^d} \inprod{\overline{V}, U} \D x 
    &= \int\limits_{\R^d} \inprod{U, W}\D x + \int\limits_{\R^d} \deriv{W} \det\left(\cnst{D} \diff\right) \trace\left([\cnst{D}\diff]^{-1} [\cnst{D} U]\right) \D x \\
    &= \int\limits_{\R^d} \inprod{U, W - \cdiv \left(\deriv{W} \det(\cnst{D}\diff) \left[\cnst{D} \diff\right]^{-1} \right)} \D x,
\end{align}
resulting in 
\begin{equation}
    [\cnst{D} \emb(\diff)]^{\ad} (W, \deriv{W}) = W - \cdiv\left(\deriv{W} \det\left(\cnst{D} \diff\right)\left[\cnst{D} \diff\right]^{-1}\right).
\end{equation} 
This allows computing the second derivative $\left[\cnst{D}\emb\left(\diff\right)\right]\left(U, V\right)$ of $\emb$ as 
\begin{align}
    \left[\cnst{D}^2\emb\left(\diff\right)\right]\left(U, V\right) 
    =& \frac{\D}{\D \epsilon} \bigg|_{\eps = 0} \Big[\cnst{D} \emb\left(\diff + \epsilon V\right) \Big]U \\
    =\quad & \left(0, \det(\cnst{D} \diff) \trace\left([\cnst{D}\diff]^{-1} [\cnst{D} U]\right)\trace\left([\cnst{D}\diff]^{-1} [\cnst{D} V]\right)\right)\\
    - &\left(0, \det(\cnst{D} \diff) \trace\left([\cnst{D}\diff]^{-1} [\cnst{D} U][\cnst{D}\diff]^{-1} [\cnst{D} V]\right)\right).
\end{align}
We now derive the dual equations of inertial motion $\ddot{\diff}_t + \Gamma^{\bpot, \emb}_{\diff_t}(\dot{\diff}_t, \dot{\diff}_t) = 0$ in the form 
\begin{align}
    [\cnst{D} \xi&\left(\diff\right)]^{\ad}\left[\cnst{D}^2\bpot \left(\xi\left(\diff\right)\right)\right][\cnst{D}\xi\left(\diff\right)] \ddot{\diff} = [\cnst{D} \xi\left(\diff\right)]^{\ad}\left[\cnst{D}^2\bpot \left(\xi\left(\diff\right)\right)\right]\\
    &\left(\left[\cnst{D}^2 \xi\left(\diff\right) \right]\left(\dot{\diff}, \dot{\diff}\right) + \left[\cnst{D}^2\bpot\left(\xi\left(\diff\right)\right)\right]^{-1}\left[\cnst{D}^3\bpot\left(\xi\left(\diff\right)\right)\right]\left(\left[\cnst{D}\xi\left(\diff\right)\right] \dot{\diff}, \left[\cnst{D}\xi\left(\diff\right)\right] \dot{\diff}\right)\right).
\end{align}
We stop tracking the $t$ dependence explicitly to simplify notation.  We first compute 
\begin{align}
   &[\cnst{D} \xi\left(\diff\right)]^{\ad}\left[\cnst{D}^2\bpot \left(\xi\left(\diff\right)\right)\right][\cnst{D}\xi\left(\diff\right)] \ddot{\diff} \\
   =&[\cnst{D} \xi\left(\diff\right)]^{\ad}\left[\cnst{D}^2\bpot \left(\xi\left(\diff\right)\right)\right]\left(\ddot{\diff}, \det\left(\cnst{D} \diff\right) \trace\left([\cnst{D}\diff]^{-1} [\cnst{D} \ddot{\diff}]\right) \right)\\ 
   =&[\cnst{D} \xi\left(\diff\right)]^{\ad}\Big( \ddot{\diff}, \alpha \trace\left([\cnst{D}\diff]^{-1} [\cnst{D} \ddot{\diff}] \right) / \det\left(\cnst{D} \diff \right) \Big) \\
   =&\ddot{\diff} - \alpha \cdiv\left( \trace\left([\cnst{D}\diff]^{-1} [\cnst{D} \ddot{\diff}]\right) \left[\cnst{D}\diff\right]^{-1}\right). 
\end{align}
Here, the operator acting on $\ddot{\diff}$ is self-adjoint and strictly positive definite, since 
\begin{equation}
        \int \limits_{\R^d} \! \inprod{- \cdiv\left( \trace\left([\cnst{D}\diff]^{-1} [\cnst{D} U]\right) \left[\cnst{D}\diff\right]^{-1}\right), V} \! \D x
        \!=\! \!
        \int \limits_{\R^d} \trace\left([\cnst{D}\diff]^{-1} [\cnst{D} U]\right) \trace\left([\cnst{D}\diff]^{-1} [\cnst{D} V]\right) \! \D x.
\end{equation}
We next compute 
\begin{align}
&\left(\left[\cnst{D}^2 \xi\left(\diff\right) \right]\left(\dot{\diff}, \dot{\diff}\right) + \left[\cnst{D}^2\bpot\left(\xi\left(\diff\right)\right)\right]^{-1}\left[\cnst{D}^3\bpot\left(\xi\left(\diff\right)\right)\right]\left(\left[\cnst{D}\xi\left(\diff\right)\right] \dot{\diff}, \left[\cnst{D}\xi\left(\diff\right)\right] \dot{\diff}\right)\right)\\
= \phantom{-} &\left(0, \det(\cnst{D} \diff) \trace\left([\cnst{D}\diff]^{-1} [\cnst{D} \dot{\diff}]\right)\trace\left([\cnst{D}\diff]^{-1} [\cnst{D} \dot{\diff}]\right)\right)\\
- &\left(0, \det(\cnst{D} \diff) \trace\left([\cnst{D}\diff]^{-1} [\cnst{D} \dot{\diff}][\cnst{D}\diff]^{-1} [\cnst{D} \dot{\diff}]\right)\right)\\
- &\left(0, 2 \det\left(\cnst{D} \diff\right) \trace\left([\cnst{D}\diff]^{-1} [\cnst{D} \dot{\diff}]\right) \trace\left([\cnst{D}\diff]^{-1} [\cnst{D} \dot{\diff}]\right) \right)\\
= - &\left(0, \det(\cnst{D} \diff) \trace\left([\cnst{D}\diff]^{-1} [\cnst{D} \dot{\diff}][\cnst{D}\diff]^{-1} [\cnst{D} \dot{\diff}]\right)\right)\\
- &\left(0, \det\left(\cnst{D} \diff\right) \trace\left([\cnst{D}\diff]^{-1} [\cnst{D} \dot{\diff}]\right) \trace\left([\cnst{D}\diff]^{-1} [\cnst{D} \dot{\diff}]\right) \right).
\end{align}
Plugging this into the dual equations of inertial motion yields, for $\trace^2\left(\cdot\right)$ the squared trace,
\begin{align}
\ddot{\diff} - \alpha \cdiv\left( \trace\left([\cnst{D}\diff]^{-1} [\cnst{D} \ddot{\diff}]\right) \left[\cnst{D}\diff\right]^{-1}\right) 
= &- \alpha \cdiv  \left(\trace\left(\left([\cnst{D}\diff]^{-1} [\cnst{D} \dot{\diff}]\right)^2\right)[\cnst{D}\diff]^{-1}\right) \\
&- \alpha \cdiv\left(\trace^2\left([\cnst{D}\diff]^{-1} [\cnst{D} \dot{\diff}]\right) [\cnst{D}\diff]^{-1}\right).
\end{align}

\subsection{In Eulerian coordinates}
We now write the Eulerian velocity field $\vct{u}$ as 
\begin{equation}
     \dot{\diff}(x) \eqqcolon \vct{u}\left(\diff(x)\right) \quad \Rightarrow \quad  
     \ddot{\diff}(x) = \left[\cnst{D} \vct{u} (\diff(x))\right] \vct{u}\left(\diff\left(x\right)\right)  + \partial_t \vct{u}(\diff(x))
\end{equation}
and the Eulerian density $\rho$ as $\rho(\diff(x)) = \left(\det \left(\cnst{D}\diff(x)\right)\right)^{-1}$.
We integrate both sides of the dual equations of motion against a vector field $V$, obtaining
\begin{align}
    & \int \limits_{\R^d} \inprod{V, \left(\ddot{\diff} - \alpha \cdiv \left(\trace\left([\cnst{D}\diff]^{-1} [\cnst{D} \ddot{\diff}]\right) [\cnst{D}\diff]^{-1}  \right)\right)} \D x\\
    & =  - \alpha \int \limits_{\R^d} \inprod{V, \cdiv\left(\trace\left(\left([\cnst{D}\diff]^{-1} [\cnst{D} \dot{\diff}]\right)^2 \right) [\cnst{D}\diff]^{-1}+\trace^2\left([\cnst{D} \dot{\diff}] [\cnst{D}\diff]^{-1}\right)[\cnst{D}\diff]^{-1} \right)} \D x
\end{align}
Integrating by parts, we obtain
\begin{align}
    &\int \limits_{\R^d} \inprod{V, \ddot{\diff}} + \alpha \trace\left([\cnst{D} V][\cnst{D}\diff]^{-1}\right) \trace\left( [\cnst{D}\diff]^{-1} [\cnst{D} \ddot{\diff}]\right)  \D x\\ 
    &= \alpha \int \limits_{\R^d} \trace\left([\cnst{D}V] [\cnst{D}\diff]^{-1}\right) \left( \trace\left(\left([\cnst{D}\diff]^{-1} [\cnst{D} \dot{\diff}]\right)^2 \right)  
    +  \trace^2\left([\cnst{D} \dot{\diff}] [\cnst{D}\diff]^{-1}\right)\right) \D x.
\end{align}
We now plug in the definition of $\vct{u}$, replacing $\dot{\diff}$ and $\ddot{\diff}$ to obtain
\begin{align}
    & \int \limits_{\R^d} \inprod{V, \left[\left[\cnst{D} \vct{u} \right] \vct{u} 
        + \partial_t \vct{u}\right]\circ{\diff}} 
        + \alpha 
            \trace \left([\cnst{D} V] [\cnst{D}\diff]^{-1}\right) 
            \trace\left( \left[\cnst{D} \left(\left[\cnst{D} \vct{u} \right] \vct{u} + \partial_t \vct{u}\right)\right]\circ{\diff}\right) 
            \D x\\ 
    &= \alpha \int \limits_{\R^d} \trace\left([\cnst{D}V] [\cnst{D}\diff]^{-1}\right) \left( \trace\left([\cnst{D} \vct{u}]^2 \circ \diff \right)  
    +  \trace^2\left([\cnst{D} \vct{u}] \circ \diff\right)\right) \D x.
\end{align}
We now write $V(x) = \tilde{V}(\diff(x))$ to obtain 
\begin{align}
    &\int \limits_{\R^d} \inprod{\tilde{V} \circ \diff, \left[\left[\cnst{D} \vct{u} \right] \vct{u} + \partial_t \vct{u}\right]\circ{\diff}} + \alpha \trace \left([\cnst{D} \tilde{V}] \circ \diff  \right) \trace\left( \left[\cnst{D} \left(\left[\cnst{D} \vct{u} \right] \vct{u} + \partial_t \vct{u}\right)\right]\circ{\diff}\right)  \D x\\ 
    &= \alpha \int \limits_{\R^d} \trace\left([\cnst{D}\tilde{V}] \circ \diff\right) \left( \trace\left([\cnst{D} \vct{u}]^2 \circ \diff \right)  
    +  \trace^2\left([\cnst{D} \vct{u}] \circ \diff\right)\right) \D x.
\end{align}
We now perform a transformation of variables with $\diff^{-1}$, resulting in multiplication with its Jacobian determinant $\det\left([\cnst{D} \diff^{-1}]\right) = \det([\cnst{D} \diff])^{-1} \circ \diff^{-1}(x) = \rho(x)$, 
\begin{align}
    & \int \limits_{\R^d} \inprod{\tilde{V}, \left[\left[\cnst{D} \vct{u} \right] \vct{u} + \partial_t \vct{u}\right] \rho} 
    + \alpha \trace \left([\cnst{D} \tilde{V}]  \right) \trace\left( \left[\cnst{D} \left(\left[\cnst{D} \vct{u} \right] \vct{u} + \partial_t \vct{u}\right)\right]\right)  \rho \D x\\ 
    =  \alpha &\int \limits_{\R^d} \trace\left([\cnst{D} \tilde{V}] \right) \left(\trace\left([\cnst{D} \vct{u}]^2 \right)   
    + \trace^2 \left([\cnst{D} \vct{u}]\right)\right) \rho \D x.
\end{align}
We again integrate by parts, obtaining
\begin{align}
    & \int \limits_{\R^d} \inprod{\tilde{V}, \left[\left[\cnst{D} \vct{u} \right] \vct{u} + \partial_t \vct{u}\right] \rho - \alpha \cdiv \left( \trace \left(\cnst{D} \left(\left[\left[\cnst{D} \vct{u} \right] \vct{u} + \partial_t \vct{u}\right]\right)\right) \rho\Id \right)} \\
    & = - \alpha \int \limits_{\R^d} \inprod{\tilde{V}, \cdiv\left(\trace\left([\cnst{D} \vct{u}]^2 \right)\rho \Id  + \trace^2\left([\cnst{D} \vct{u}] \right) \rho \Id \right)} \D x.
\end{align}
The general choice of $\tilde{V}$ and the fundamental lemma of the calculus of variations yield
\begin{equation}
    \left( (\cdot) \rho - \alpha \cdiv\left( \rho \trace\left(\cnst{D} \left(\cdot \right)\right) \Id \right)\right) \left(\left[\cnst{D} \vct{u} \right] \vct{u} + \partial_t \vct{u}\right)
    = -\alpha \cdiv\left(\trace\left([\cnst{D}\vct{u}]^2\right) \rho \Id  
      + \trace^2\left([\cnst{D}\vct{u}]\right)  \rho \Id\right).
\end{equation}
All matrix divergences above are taken over multiples of the identity matrix.
Thus, we can safely replace them with row-wise divergences, obtaining
\begin{equation}
     \left(\left[\cnst{D} \vct{u} \right] \vct{u} + \partial_t \vct{u}\right)
    =- \alpha \left( (\cdot) \rho - \alpha \rdiv\left( \rho \trace\left(\cnst{D} \left(\cdot \right)\right) \Id \right)\right)^{-1}  
        \rdiv\left(\trace\left([\cnst{D}\vct{u}]^2\right) \rho \Id + \trace^2\left([\cnst{D}\vct{u}]\right) \rho \Id\right).
\end{equation}
We multiply both sides with $\rho$ and apply the conservation of mass to obtain
\begin{align}
     \partial_t \left(\rho \vct{u}\right) + \rdiv\left(\rho \vct{u} \otimes \vct{u}\right) 
    = - \alpha \rho \Big( (\cdot) \rho -& \alpha \rdiv\left( \rho \trace\left(\cnst{D} \left(\cdot \right)\right) \Id \right)\Big)^{-1}\\  
        &\rdiv\left(\trace\left([\cnst{D}\vct{u}]^2 \right) \rho \Id + \trace^2\left([\cnst{D}\vct{u}]\right)  \rho \Id\right).
\end{align}
\subsubsection{In conservation form} We now want to bring the regularized equations into conservation form. 
As in \cref{sec:conservation_1d}, we pull divergences and multiplications with $\rho$ through the inverse to obtain
\begin{align}
&- \alpha \rho \left( (\cdot) \rho - \alpha \rdiv\left( \rho \trace\left(\cnst{D} \left(\cdot \right)\right) \Id \right)\right)^{-1}
        \rdiv\left(\trace\left([\cnst{D}\vct{u}]^2 \right) \rho \Id + \trace^2\left([\cnst{D}\vct{u}]\right)  \rho \Id\right)\\
=
&- \alpha \left( (\cdot) - \alpha \rdiv \left(\rho \trace\left(\cnst{D} \left(\rho^{-1}\left(\cdot \right)\right)\right) \Id \right)\right)^{-1}
\rdiv \left(\trace\left([\cnst{D}\vct{u}]^2\right) \rho \Id 
            + \trace^2\left([\cnst{D}\vct{u}]\right) \rho \Id\right)\\
= 
&- \alpha \left( \rdiv^{-1}(\cdot) - \alpha \rho \trace\left(\cnst{D} \left(\rho^{-1}\left(\cdot \right)\right)\right) \Id \right)^{-1}
    \left(\trace\left([\cnst{D}\vct{u}]^2\right) \rho \Id 
        + \trace^2\left([\cnst{D}\vct{u}]\right) \rho \Id\right)\\
=
&- \alpha \rdiv\left( (\cdot) - \alpha \rho \trace\left(\cnst{D} \left(\rho^{-1}\rdiv\left(\cdot \right)\right)\right) \Id \right)^{-1}
    \left(\trace\left([\cnst{D}\vct{u}]^2\right) \rho \Id 
        + \trace^2\left([\cnst{D}\vct{u}]\right) \rho \Id\right)\\
=
&- \alpha \rdiv\left( \rho^{-1} (\cdot) - \alpha \trace\left(\cnst{D} \left(\rho^{-1}\rdiv\left(\cdot \right)\right)\right) \Id \right)^{-1}
    \left(\trace\left([\cnst{D}\vct{u}]^2\right) \Id 
        + \trace^2\left([\cnst{D}\vct{u}]\right) \Id\right).
\end{align}
In other words, the additional momentum flux satisfies 
\begin{equation}
    \rho^{-1} \mtx{F} - \alpha \trace\left(\cnst{D} \left(\rho^{-1} \rdiv \mtx{F}\right) \right)\mtx{I} 
    = \alpha \left(\trace\left([\cnst{D}\vct{u}]^2 \right) \Id 
        + \trace^2\left([\cnst{D}\vct{u}]\right) \Id\right).
\end{equation}
It follows that $\mtx{F}$ is a multiple of the identity of the form $\mtx{F} = \Sigma \Id$, resulting in
\begin{equation}
    \rho^{-1} \Sigma - \alpha \trace\left(\cnst{D} \left(\rho^{-1} \nabla \Sigma \right) \right)
    = \alpha \left(\trace\left([\cnst{D}\vct{u}]^2 \right) 
        + \trace^2\left([\cnst{D}\vct{u}]\right)\right).
\end{equation}
The trace of the Jacobian is just the (vector) divergence, resulting in 
\begin{equation}
    \rho^{-1} \Sigma - \alpha \divergence \left(\rho^{-1} \nabla \Sigma \right) 
    = \alpha \left(\trace\left([\cnst{D}\vct{u}]^2 \right) 
        + \trace^2\left([\cnst{D}\vct{u}]\right)\right).
\end{equation}
After adding the pressure and external forces, this yields the information geometric regularization of the multidimensional barotropic Euler equation 
\begin{equation}
    \label{eqn:infgeo_reg_flux_dd}
     \begin{cases}
          \partial_t
          \begin{pmatrix}
               \rho \vct{u} \\
               \rho
          \end{pmatrix}
          + \rdiv 
          \begin{pmatrix}
               \rho \vct{u} \otimes \vct{u} + \left(P(\rho) + \Sigma\right) \Id\\ 
               \rho \vct{u}
          \end{pmatrix}
          = 
          \begin{pmatrix} 
               \vct{f}\\
               0
          \end{pmatrix} \\
          \rho^{-1} \Sigma - \alpha \divergence(\rho^{-1} \nabla \Sigma) = \alpha  \left(\trace^2\left([\cnst{D} \vct{u}]\right) + \trace\left([\cnst{D} \vct{u}]^2\right) \right).
     \end{cases}
\end{equation}
To rigorously justify the above manipulations, it is easiest to first establish the existence of a solution $\Sigma$ of the scalar elliptic PDE in \cref{eqn:infgeo_reg_flux_dd} and then trace the arguments backwards to show that $\vct{F} = \Sigma \Id$ satisfies the equations above.

\section{Numerical experiments}
\label{sec:numerics}
\begin{figure}
     \begin{tikzpicture}
\begin{groupplot}[
	group style={group size=2 by 1,
	horizontal sep=0.075\textwidth,
    vertical sep=0.5cm,},
	]
	\nextgroupplot[
		title={$t = 0.0875$},
		standard,
		xlabel={$x$},
		ylabel={$u(x, t)$},
		height=0.415\textwidth,
		width=0.415\textwidth,
        ymin=-3.0,
        ymax=3.0,
        xmin=0.4,
        xmax=0.6,
		xtick={0.45,0.5, 0.55},
		ytick={-2, -1, 0, 1, 2},
		enlarge x limits=0.,
		enlarge y limits=0.,
		legend image post style={scale=1.5},
		legend style={
			at={(1.10, 1.25)},inner sep=3pt,anchor=north,legend columns= -1,cells ={align = right}, draw=none,fill=none}, 
		]	

		\addlegendimage{steelblue, line width=1.8pt}
		\addlegendimage{orange, line width=1.8pt}
		\addlegendimage{joshua, line width=1.8 pt}
		\addlegendimage{rust, line width=1.8 pt}

		\legend{Reference solution, LF, LW, LW + IGR};

		\addplot[rust, line width=1.8pt] table[x index={2},y index={0}] {csv/1d_sine/igr_T1.csv};
		\addplot[orange, line width=1.8pt] table[x index={2},y index={0}] {csv/1d_sine/lf_T1.csv};
		\addplot[joshua, line width=1.8pt] table[x index={2},y index={0}] {csv/1d_sine/rlw_T1.csv};
		\addplot[steelblue, line width=1.8pt] table[x index={2},y index={0}] {csv/1d_sine/ref_T1.csv};

		\nextgroupplot[
			standard,
			title={$t = 0.75$},
			xlabel={$x$},
			height=0.415\textwidth,
			width=0.415\textwidth,
			ymin=-1.0,
			ymax=1.0,
			xmin=0.4,
			xmax=0.6,
			xtick={0.45,0.5, 0.55},
			ytick={-1.0, -0.5, 0.0, 0.5, 1.0},
			enlarge x limits=0.,
			enlarge y limits=0.,
			legend image post style={scale=1.5},
			legend style={
				at={(2.15,1.125)},inner sep=4.5pt,anchor=south,legend columns=9,legend cell align={left}, draw=none,fill=none}
			]	
	
			\addplot[rust, line width=1.8pt] table[x index={2},y index={0}] {csv/1d_sine/igr_T2.csv};
			\addplot[orange, line width=1.8pt] table[x index={2},y index={0}] {csv/1d_sine/lf_T2.csv};
			\addplot[joshua, line width=1.8pt] table[x index={2},y index={0}] {csv/1d_sine/rlw_T2.csv};
			\addplot[steelblue, line width=1.8pt] table[x index={2},y index={0}] {csv/1d_sine/ref_T2.csv};
\end{groupplot}
\end{tikzpicture}
     \vspace{-8mm}
     \caption{\textbf{IGR on Euler:} Lax-Wendroff (LW) with information geometric regularization (IGR) avoids oscillations in LW and the oversmoothing of Lax-Friedrichs (LF).}
     \label{fig:1d_sine}
\end{figure}

\begin{figure}
     \vspace{-4mm}
     \begin{tikzpicture}
\begin{groupplot}[
	group style={group size=2 by 1,
	horizontal sep=0.075\textwidth,
    vertical sep=0.5cm,},
	]
			\nextgroupplot[
				standard,
				xlabel={$x$},
				ylabel={$\rho(x, 2)$},
				height=0.415\textwidth,
				width=0.415\textwidth,
				ymin=0.5,
				ymax=4.5,
        		xmin=3.0,
        		xmax=9.0,
				xtick={4.0, 5.0, 6.0, 7.0, 8.0},
				ytick={1.0, 2.0, 3.0, 4.0},
				yticklabels={{$\phantom{-}1$, $2$, $3$, $4$}},
				enlarge x limits=0.,
				enlarge y limits=0.,
				legend image post style={scale=1.5},
				legend style={
					at={(0.95, 1.1)},inner sep=3pt,anchor=north,legend columns= -1,cells ={align = right}, draw=none,fill=none}, 
				]	

		\addlegendimage{steelblue, line width=1.8pt}
		\addlegendimage{orange, line width=1.8pt}
		\addlegendimage{joshua, line width=1.8 pt}
		\addlegendimage{rust, line width=1.8 pt}


			\addplot[steelblue, line width=1.8pt] table[x index={2},y index={1}] {csv/shu_osher/ref.csv};
			\addplot[orange, line width=1.8pt] table[x index={2},y index={1}] {csv/shu_osher/lf.csv};
			\addplot[joshua, line width=1.8pt] table[x index={2},y index={1}] {csv/shu_osher/rlw.csv};
			\addplot[rust, line width=1.8pt] table[x index={2},y index={1}] {csv/shu_osher/igr.csv};
		
	\nextgroupplot[
		standard,
 		xlabel={$x$},
		height=0.415\textwidth,
		width=0.415\textwidth,
        ymin=3.5,
        ymax=4.1,
        xmin=3.75,
        xmax=7.75,
		xtick={4.0, 5.0, 6.0, 7.0},
		ytick={3.6, 3.8, 4.0},
		enlarge x limits=0.,
		enlarge y limits=0.,
	]
		\addplot[steelblue, line width=1.8pt] table[x index={2},y index={1}] {csv/shu_osher/ref.csv};
		\addplot[orange, line width=1.8pt] table[x index={2},y index={1}] {csv/shu_osher/lf.csv};
		\addplot[joshua, line width=1.8pt] table[x index={2},y index={1}] {csv/shu_osher/rlw.csv};
		\addplot[rust, line width=1.8pt] table[x index={2},y index={1}] {csv/shu_osher/igr.csv};
	\end{groupplot}
\end{tikzpicture}
     \vspace{-8mm}
     \caption{\textbf{Interaction of shock and sound waves:} Information geometric regularization (IGR) with Lax-Wendroff (LW) correctly predicts the magnitude of the sound waves. Plain LW overestimates them while Lax-Friedrichs (LF) underestimates them.}
     \label{fig:shu_osher}
     \vspace{-4mm}
\end{figure}

\subsection{Overview} This section illustrates the behavior of solutions to the modified equations and provides proof of the concept of using them for numerical purposes. 
We emphasize that the results in this section are not yet fair comparisons to the state-of-the-art in simulating compressible flows.
Such would require stronger baselines and the consideration of practical aspects like computational cost, memory, and robustness. 
This is beyond the scope of this work and will be addressed in the future.
We focus on Lax-Friedrichs and Lax-Wendroff methods as classical examples of the tradeoff between dissipation and Gibbs oscillations that is at the heart of all shock capturing. 
We use explicit methods and do not observe any additional CFL restrictions due to IGR.
The code used to generate this section's results is available at \url{https://github.com/f-t-s/information_geometric_regularization_of_barotropic_euler}.

\begin{remark}[Boundary conditions]
     We derive IGR in $\R^d$ but practical computation requires imposing boundary conditions. 
     In the numerical experiments of this section, we use periodic boundary conditions that are implemented by periodic extension (instead of the value at $x_{N + k}$ we use that at $x_k$).
     Preliminary results suggest that our regularization also applies to other cases like Dirichlet and Neumann conditions.
     This requires Neumann boundary conditions for the elliptic problem defining $\Sigma$. 
\end{remark}

\subsection{The unidimensional case}
We begin with an application to the unidimensional Euler equation. 
We use the pressure $P(\rho) = \rho^{1.4}$ and compare the solutions obtained by the Lax-Friedrichs (LF) scheme, the Richtmyer Lax-Wendroff (LW) scheme \cite{leveque1992numerical}, and LW applied to \cref{eqn:infgeo_reg_flux_1d} with $\alpha = 20 (\Delta x)^2$ (denoted as LW + IGR).
In terms of the function $\vct{F}$ mapping the state vector $\vct{q}$ to the flux,  the LF scheme is
\begin{equation}
\vct{q}_{i}^{k + 1} = \frac{1}{2}\left(\vct{q}_{i + 1}^{k} + \vct{q}_{i - 1}^{k}\right) - \frac{\Delta t}{2 \Delta x}\left(\vct{F}\left(\vct{q}_{i + 1}^{k}\right) - \vct{F}\left(\vct{q}_{i-1}^{k}\right)\right)
\end{equation}
and the LW scheme has the form 
\begin{equation}
\begin{split}
\vct{q}_{i + 1 / 2}^{k + 1/2} &= \frac{1}{2}\left(\vct{q}_{i + 1}^{k} + \vct{q}_{i}^{k}\right) - \frac{\Delta t}{2 \Delta x}\left(\vct{F}\left(\vct{q}_{i + 1}^{k}\right) - \vct{F}\left(\vct{q}_{i}^{k}\right)\right)\\
\vct{q}_{i}^{k + 1} &= \vct{q}_{i}^k - \frac{\Delta t}{\Delta x} \left(\vct{F}\left(\vct{q}^{k + 1 / 2}_{i + 1/2}\right) -  \vct{F}\left(\vct{q}^{k + 1 / 2}_{i - 1/2}\right)\right).
\end{split}
\end{equation}
Evaluating the function $\vct{F}$ in the case \cref{eqn:reg_euler_intro} requires computing $\Sigma$. 
We discretize the elliptic problem defining it using a tridiagonal Laplacian stencil, compute the right-hand side using central differences, and use sparse Gaussian elimination to solve for $\Sigma$, at cost linear in the size of the mesh.
As is well known, the viscous regularization of LF causes excessive smoothness and the dispersive errors of LW cause spurious ``Gibbs'' oscillations near shocks. 
We are testing the hypothesis that IGR avoids these Gibbs oscillations by ensuring that the underlying continuous problem has a regular solution, without adding excessive smoothing. 
IGR tames shock formation and insures that the regularity of the initial condition is maintained, it does not aid in problems with irregular initial conditions, like the advection of a step function. 
Initial conditions with step functions can be handled by smearing the singularity over a few grid points. 
For the sake of comparability across grid sizes, we do so by replacing jumps in the initial condition with sigmoidal blending or smooth cutoff functions. 
In practical applications, it is easier to instead apply a few iterations of heat kernel smoothing to the initial and boundary conditions.

\subsubsection{Sine wave initial condition}
\label{sec:sine_wave_1d}
For our first experiment, we choose initial conditions $u(x, 0) = 3 \sin(2 \pi x), \rho(x, 0) \equiv 1$, discretize the domain $[0, 1]$ with 500 grid points and use periodic boundary conditions. 
As shown in \cref{fig:1d_sine}, LW + IGR avoids both the excessive dissipation of LF and the Gibbs oscillations of LW.
We compare to a reference solution obtained by Lax-Friedrichs with $2 \times 10^4$ grid points.

\begin{figure}
     \begin{tikzpicture}
\begin{groupplot}[
	group style={group size=4 by 1,
	horizontal sep=0.0075\textwidth,
    vertical sep=0.5cm,},
	]
	\nextgroupplot[
		title={$t = 0.055$},
		standard,
		xlabel={$x$},
		ylabel={$u(x, t)$},
		height=0.19\textwidth,
		width=0.19\textwidth,
        ymin=-3.0,
        ymax=3.0,
        xmin=0.0,
        xmax=1.0,
		xtick={0.25,0.5, 0.75},
		ytick={-2, -1, 0, 1, 2},
		enlarge x limits=0.,
		enlarge y limits=0.,
		legend image post style={scale=1.5},
		legend style={
			at={(2.25, 1.70)},inner sep=3pt, anchor=north,legend columns= 4,cells ={align = right}, draw=none,fill=none, /tikz/every even column/.append style={column sep=0.025\textwidth}}, 
		]	

		\addlegendimage{darksky, line width=1.8pt}
		\addlegendimage{seagreen, mark=triangle, only marks, mark size=2.5pt, thick}
		\addlegendimage{seagreen, mark=o, only marks, mark size=2.5pt, thick}
		\addlegendimage{steelblue, line width=1.8 pt}
		\addlegendimage{rust, line width=1.8pt}
		\addlegendimage{orange, mark=triangle, only marks, mark size=2.5pt, thick}
		\addlegendimage{orange, mark=o, only marks, mark size=2.5pt, thick}
		\addlegendimage{joshua, line width=1.8 pt}

		\legend{Smooth Reference, Smooth $u$, Smooth $\rho$, Slope $-2$, Shocked Reference, Shocked $u$, Shocked $\rho$, Slope $-1$};

		\addplot[darksky, line width=1.8pt] table[x index={0},y index={2}] {csv/convergence_study/smooth_ref_sol.csv};

		\nextgroupplot[
			standard,
			title={$t = 0.1$},
			xlabel={$x$},
			height=0.19\textwidth,
			width=0.19\textwidth,
			ymin=-3.0,
			ymax=3.0,
			xmin=0.0,
			xmax=1.0,
			yticklabels=\empty,
			xtick={0.25,0.5, 0.75},
			ytick={-2, -1, 0, 1, 2},
			enlarge x limits=0.,
			enlarge y limits=0.,
			]	
	
			\addplot[rust, line width=1.8pt] table[x index={0},y index={2}] {csv/convergence_study/shock_ref_sol.csv};
			
		\nextgroupplot[
			xshift=0.1\textwidth,
			standard,
			title={$t = 0.055$},
			xlabel={$- \log_{10}\big(\sqrt{\alpha}\big)$},
			ylabel={$\log_{10}\left(\text{relative $\|\cdot \|_1$ error}\right)$},
			height=0.19\textwidth,
			width=0.19\textwidth,
			ymin=-6.0,
			ymax=1.0,
			xmin=0.5,
			xmax=4.5,
			xtick={2.0, 4.0, 6.0, 8.0, 1.0},
			ytick={-1, -3.0, -5.0},
			enlarge x limits=0.,
			enlarge y limits=0.,
			]	
	
			\addplot[steelblue, line width=1.8pt, samples=100, domain=0:10] {-2 * x + 2.5};
			\addplot[seagreen, mark = o, only marks, mark size=2.5pt, thick] table[x expr={-log10(sqrt(\thisrowno{0}))},y expr={log10(\thisrowno{3})}] {csv/convergence_study/smooth_errors.csv};
			\addplot[seagreen, mark = triangle, only marks, mark size=2.5pt, thick] table[x expr={-log10(sqrt(\thisrowno{0}))},y expr={log10(\thisrowno{1})}] {csv/convergence_study/smooth_errors.csv};

		\nextgroupplot[
			standard,
			title={$t = 0.1$},
			xlabel={$- \log_{10}\big(\sqrt{\alpha}\big)$},
			height=0.19\textwidth,
			width=0.19\textwidth,
			ymin=-6.0,
			ymax=1.0,
			xmin=0.5,
			xmax=4.5,
			yticklabels=\empty,
			xtick={2.0, 4.0, 6.0, 8.0, 1.0},
			ytick={-1, -3.0, -5.0},
			enlarge x limits=0.,
			enlarge y limits=0.,
			]	
	
			\addplot[joshua, line width=1.8pt, samples=100, domain=0:10] {-x + 1};
			\addplot[orange, mark = o, only marks, mark size=2.5pt, thick] table[x expr={-log10(sqrt(\thisrowno{0}))},y expr={log10(\thisrowno{3})}] {csv/convergence_study/shock_errors.csv};
			\addplot[orange, mark = triangle, only marks, mark size=2.5pt, thick] table[x expr={-log10(sqrt(\thisrowno{0}))},y expr={log10(\thisrowno{1})}] {csv/convergence_study/shock_errors.csv};
\end{groupplot}
\end{tikzpicture}
     \vspace{-3mm}
     \caption{\textbf{Convergence study.} We solve the problem of \cref{sec:sine_wave_1d} with $8\cdot 10^4$ grid points, for varying $\alpha$.
     As a function of $\sqrt{\alpha}$, we observe first order convergence in the presence of shocks and second order otherwise.}
     \label{fig:convergence_study}
\end{figure}

\subsubsection{Interaction of shock and sound waves} We next consider the problem 
\begin{equation}
     u(x, 0) \approx 
     \begin{cases}
          3, &\text{ if } x < 2,\\
          0, &\text{ else,}
     \end{cases} 
     \quad 
     \rho(x, 0) \approx 
     \begin{cases}
          2, &\text{ if } x < 2\\
          1 + 0.2 \sin( 2 \pi (25 x)  / 8), &\text{ if } 2 \leq x \leq 10
     \end{cases} 
\end{equation}
that features the interaction of a right-moving shock wave and a left-moving sound wave. 
We discretize the domain $[0, 10]$ using 500 points and add boundary padding equal to $(u(\cdot,0),\rho(\cdot, 0))$. 
We use smooth cut-off functions for the transitions at $x = 2$ and at the wraparound of our periodic boundary conditions. 
We compare LF, LW, and LW + IGR to a reference solution computed with LF using $4 \times 10^4$ points on $[0, 10]$.  
As shown in \cref{fig:shu_osher}, once the two waves interact, LF damps out the sound wave while LW greatly overestimates its magnitude. 
LW + IGR estimates it correctly.
\begin{figure}
     \centering
     \includegraphics[width=\textwidth]{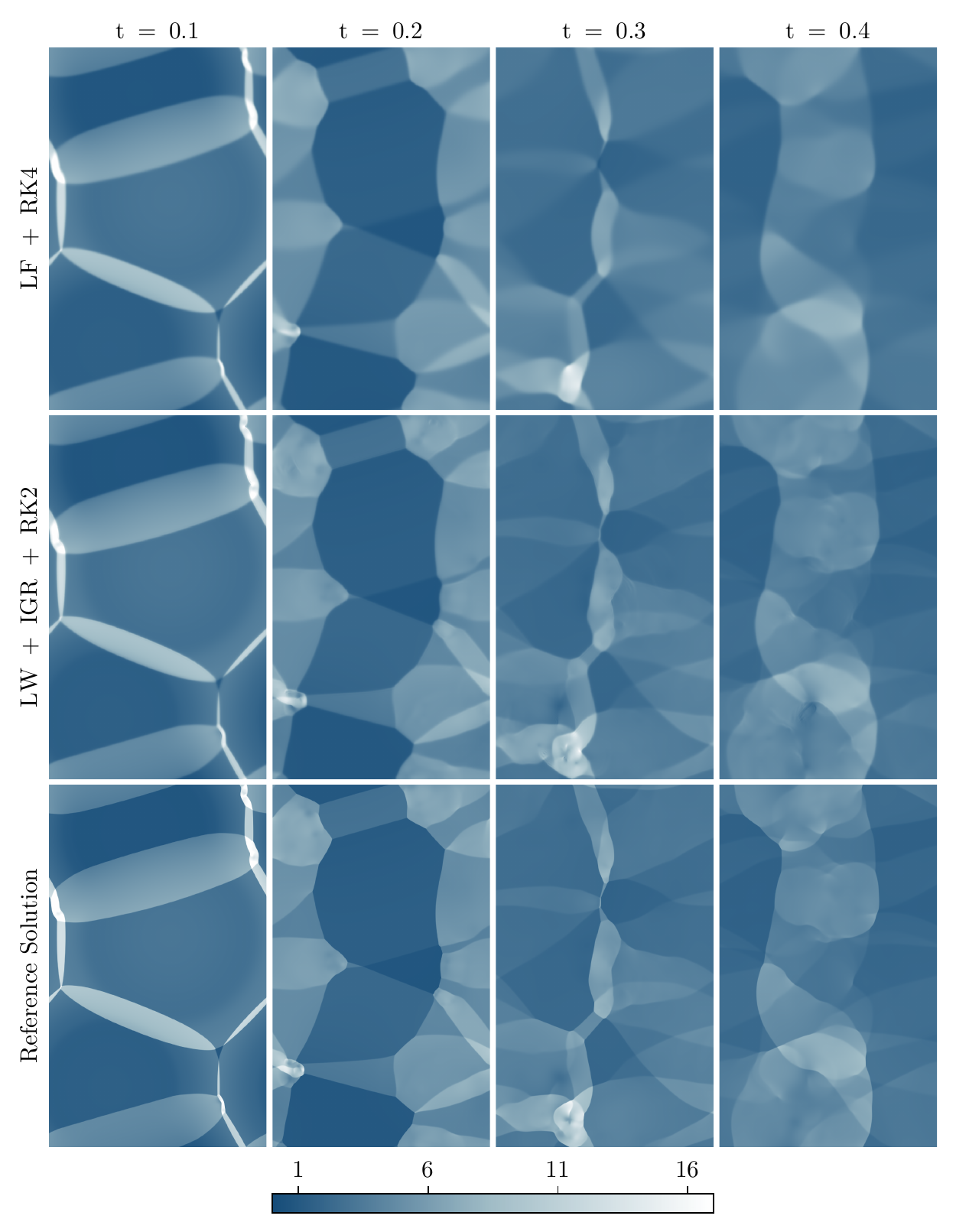}
     \vspace{-8 mm}
     \caption{\textbf{Blast waves.} Density plots of Lax-Friedrichs (LF) and Lax-Wendroff with IGR on a $432 \times 720$ grid to a reference solution computed by LF on a $4320 \times 7200$ grid.}
     \label{fig:blast_waves_density} 
\end{figure}

\subsubsection{Convergence study}
We study the convergence of highly resolved IGR solutions to those of the nominal equation, using the sine initial condition from \cref{sec:sine_wave_1d}.  
As illustrated in \cref{fig:1d_sine}, as a function of $\sqrt{\alpha}$, IGR shows first-order accuracy in the presence of shocks and second-order accuracy, otherwise.
Second order convergence in $\sqrt{\alpha}$ (and thus, in practice, in $\Delta x$) might seem disappointing given that methods of arbitrary order of convergence exist. 
But the order of accuracy of a numerical method only captures its asymptotic behavior when progressively overresolving a given solution. 
In our opinion, the concept of \emph{resolution} that measures the number of grid points required for a reasonable approximation of a given solution is the more meaningful measure of quality \cite{lax1978accuracy}.
As shown in \cref{fig:IGR_v_LAD}, IGR is much less prone to damping out oscillatory solutions than equivalent viscous regularization. 
We believe that this property will make it prove useful in simulating complex flows.

\subsection{The multidimensional case}
We now investigate the two-dimensional version of \cref{eqn:infgeo_reg_flux_dd}. 
We use the pressure $P(\rho) = \rho^{1.4}$ and multidimensional variants of the Lax-Friedrichs (LF) and Richtmyer Lax-Wendroff (LW) schemes. 
For $\vct{F}_x, \vct{F}_y$ the fluxes in the $x$ and $y$ directions, respectively, the LF scheme is given by 
\begin{equation}
\begin{split}
\vct{q}_{i,j}^{k + 1} = \frac{1}{4}\left(\vct{q}_{i + 1, j}^{k} + \vct{q}_{i - 1, j}^{k} + \vct{q}_{i, j + 1}^{k} + \vct{q}_{i, j - 1}^{k}\right) 
&- \frac{\Delta t}{2 \Delta x}\left(\vct{F}_x\left(\vct{q}_{i + 1, j}^{k}\right) - \vct{F}_x\left(\vct{q}_{i - 1, j}^{k}\right)\right) \\
&- \frac{\Delta t}{2 \Delta y}\left(\vct{F}_y\left(\vct{q}_{i, j + 1}^{k}\right) - \vct{F}_y\left(\vct{q}_{i, j - 1}^{k}\right)\right)
\end{split}
\end{equation}
and the LW scheme is given by
\begin{equation}
\begin{split}
\vct{q}_{i + 1 / 2, j + 1 / 2}^{k + 1/2} =& \frac{1}{2}\left(\vct{q}_{i + 1, j + 1}^{k} + \vct{q}_{i + 1, j}^{k} + \vct{q}_{i, j + 1}^{k} + \vct{q}_{i, j}^{k}\right)  \\
& - \frac{\Delta t}{4 \Delta x}\left(\vct{F}_x\left(\vct{q}_{i + 1, j}^{k}\right) + \vct{F}_x\left(\vct{q}_{i + 1, j + 1}^{k}\right)- \vct{F}_x\left(\vct{q}_{i, j}^{k} \right) - \vct{F}_x\left(\vct{q}_{i, j + 1}^{k}\right)\right)\\
& - \frac{\Delta t}{4 \Delta y}\left(\vct{F}_y\left(\vct{q}_{i, j + 1}^{k}\right) + \vct{F}_y\left(\vct{q}_{i + 1, j + 1}^{k}\right)- \vct{F}_y\left(\vct{q}_{i, j}^{k} \right) - \vct{F}_y\left(\vct{q}_{i + 1, j}^{k}\right)\right),
\end{split}
\end{equation}
\begin{equation}
\begin{split}
\vct{q}_{i, j}^{k + 1} = \vct{q}_{i}^k 
& - \frac{\Delta t}{4 \Delta x}\left(\vct{F}_x\left(\vct{q}_{i + 1/2, j + 1/2}^{k + 1/2}\right) + \vct{F}_x\left(\vct{q}_{i + 1/2, j}^{k + 1/2}\right)- \vct{F}_x\left(\vct{q}_{i - 1/2, j}^{k + 1/2} \right) - \vct{F}_x\left(\vct{q}_{i - 1 / 2, j - 1 / 2}^{k + 1/2}\right)\right)\\
& - \frac{\Delta t}{4 \Delta x}\left(\vct{F}_y\left(\vct{q}_{i + 1/2, j + 1/2}^{k + 1/2}\right) + \vct{F}_y\left(\vct{q}_{i, j + 1/2}^{k + 1/2}\right)- \vct{F}_y\left(\vct{q}_{i, j - 1/2}^{k + 1/2} \right) - \vct{F}_y\left(\vct{q}_{i - 1 / 2, j - 1 / 2}^{k + 1/2}\right)\right).
\end{split}
\end{equation}
\begin{figure}
     \centering
     \includegraphics[width=\textwidth]{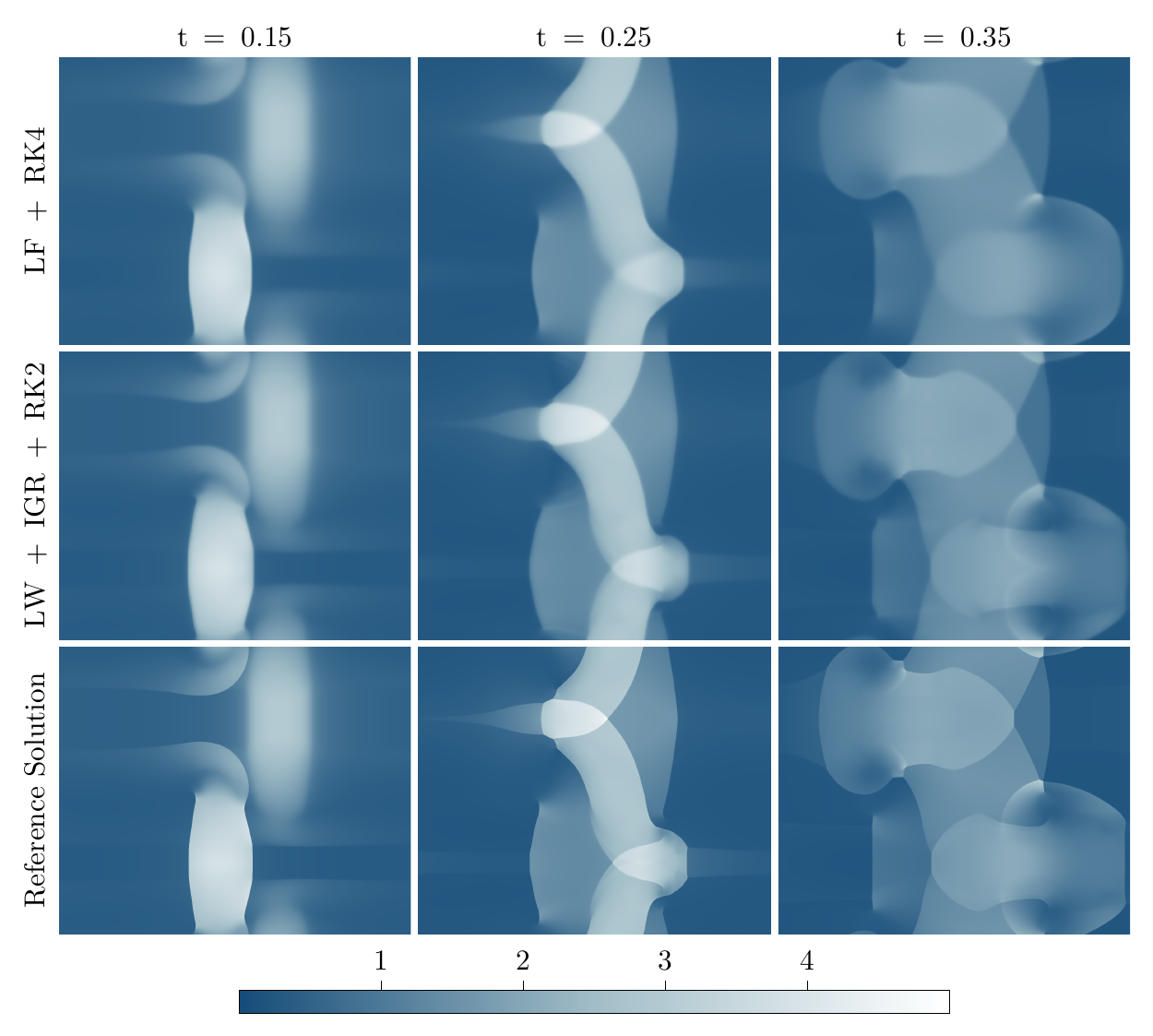}
     \vspace{-8 mm}
     \caption{\textbf{Shock collision.}\!\! We plot densities for Lax-Friedrichs and Lax-Wendroff with IGR on a $600 \times 500$ grid and a reference solution (Lax-Friedrichs, $1.2 \cdot 10^4 \times 10^4$ grid).}
     \label{fig:sine_waves_density} 
     \vspace{-4 mm}
\end{figure}
\begin{figure}
     \begin{tikzpicture}
\begin{groupplot}[
	group style={group size=3 by 1,
	horizontal sep=0.2 cm,
    vertical sep= 2 cm,},
	]
	\nextgroupplot[
		standard,
		title={Kinetic Energy},
		xlabel={$t$},
		ylabel={Energy},
		height=0.22\textwidth,
		width=0.28\textwidth,
        ymin= 0.0,
        ymax= 1.50,
        xmin=0.0,
        xmax=0.4,
		legend style={
			at={(1.50, 1.40)},inner sep=3pt,anchor=north,legend columns= -1,cells ={align = left}, draw=none,fill=none},
		]	
		\addlegendimage{steelblue, line width=2.8pt}
		\addlegendimage{orange, line width=2.8pt, mark size=2.7pt}
		\addlegendimage{rust, line width=2.8 pt, mark size=2.7, , dash pattern=on 4pt off 4pt}

        \legend{Reference solution, LF + RK4 , LW + IGR + RK2};

        \addplot[orange, line width=2.8pt] table[x index={3},y index={0}] {csv/two_sine_waves/energy_lf_two_sine_waves.csv};
        \addplot[steelblue, line width = 2.8pt] table[x index={3},y index={0}] {csv/two_sine_waves/energy_ref_two_sine_waves.csv};
        \addplot[rust, line width=2.8pt, dash pattern=on 4pt off 4pt] table[x index={3},y index={0}] {csv/two_sine_waves/energy_igr_two_sine_waves.csv};

	\nextgroupplot[
		standard,
		title={Potential Energy},
		height=0.22\textwidth,
		width=0.28\textwidth,
		ylabel={},
        ymin= 0.0,
        ymax= 1.50,
        xmin=0.0,
        xmax=0.4,
		xlabel={$t$},
		yticklabels={},
		scaled y ticks = false,
		legend style={
			at={(1.0, 1.2)},inner sep=3pt,anchor=north,legend columns= -1,cells ={align = right}, draw=none,fill=none, font = \Large},
		]	

        \addplot[orange, line width=2.8pt] table[x index={3},y index={1}] {csv/two_sine_waves/energy_lf_two_sine_waves.csv};
        \addplot[steelblue, line width = 2.8pt] table[x index={3},y index={1}] {csv/two_sine_waves/energy_ref_two_sine_waves.csv};
        \addplot[rust, line width=2.8pt, dash pattern=on 4pt off 4pt] table[x index={3},y index={1}] {csv/two_sine_waves/energy_igr_two_sine_waves.csv};

	\nextgroupplot[
		standard,
		title={Total Energy},
		height=0.22\textwidth,
		width=0.28\textwidth,
		ylabel={},
 		xlabel={$t$},
        ymin= 0.0,
        ymax= 1.50,
        xmin=0.0,
        xmax=0.4,
		yticklabels={},
		scaled y ticks = false,
		]	

        \addplot[orange, line width=2.8pt] table[x index={3},y index={2}] {csv/two_sine_waves/energy_lf_two_sine_waves.csv};
        \addplot[steelblue, line width = 2.8pt] table[x index={3},y index={2}] {csv/two_sine_waves/energy_ref_two_sine_waves.csv};
        \addplot[rust, line width=2.8pt, dash pattern=on 4pt off 4pt] table[x index={3},y index={2}] {csv/two_sine_waves/energy_igr_two_sine_waves.csv};
\end{groupplot}
\end{tikzpicture}
     \vspace{-4 mm}
     \caption{\textbf{Inviscid regularization.} We plot the energy dissipation of Lax-Friedrichs and Lax-Wendroff with IGR on a $600 \times 500$ grid, and a reference solution computed by Lax-Friedrichs on a $10^4 \times 1.2 \cdot 10^4$ grid. 
     Contrary to the Lax-Friedrichs solution, the IGR solution shows no dissipative bias.
     (Remember that $E_{\mathrm{kin}} = \frac{1}{2} \int \rho u^2 \D \vct{x}$, $E_{\mathrm{pot}} =  \int a \rho^\gamma / (\gamma - 1) \D \vct{x} - \int a / (\gamma - 1) \D \vct{x}$, and $E_{\mathrm{total}} = E_{\mathrm{kin}} + E_{\mathrm{pot}}$.)}
     \label{fig:2d_euler_sine_wave_energy}
\end{figure}
\begin{figure}
     \centering
     \includegraphics[width=\textwidth]{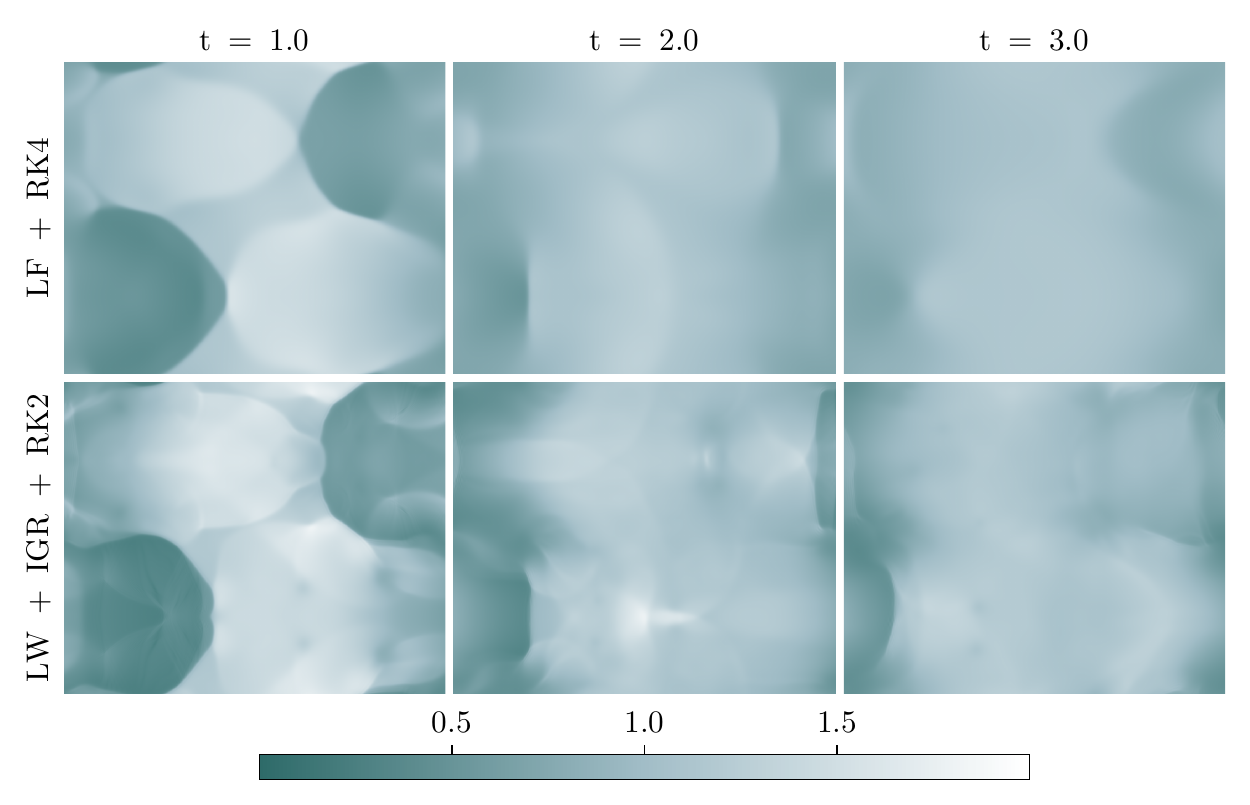}
     \vspace{-8 mm}
     \caption{\textbf{Preservation of density gradients.} When simulating the problem in \cref{fig:sine_waves_density} over long time horizons, viscous regularization eventually loses the sharpness of the density gradients. In contrast, IGR maintains the steepness of density gradients.}
     \label{fig:sine_waves_longtime_density} 
\end{figure}
\begin{figure}
     \centering
     \vspace{-2 mm}
     \includegraphics[width=\textwidth]{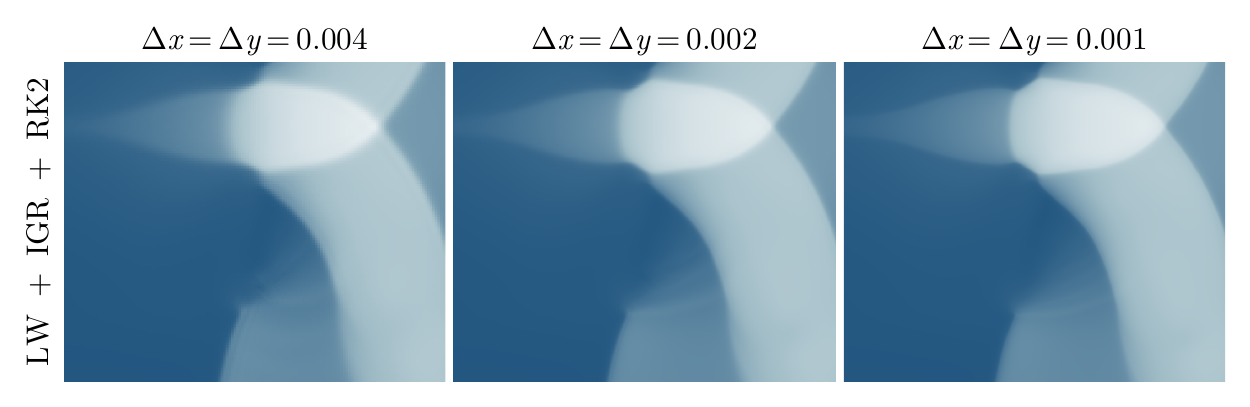}
     \vspace{-8 mm}
     \caption{\textbf{Mesh convergence.} With increasing spatial resolution and for constant $\alpha = 5 \cdot 10^{-5}$, the solutions appear to converge to a well-defined PDE solution.}
     \label{fig:sine_waves_mesh_convergence} 
     \vspace{-2 mm}
\end{figure}
\begin{figure}
     \centering
     \includegraphics[width=\textwidth]{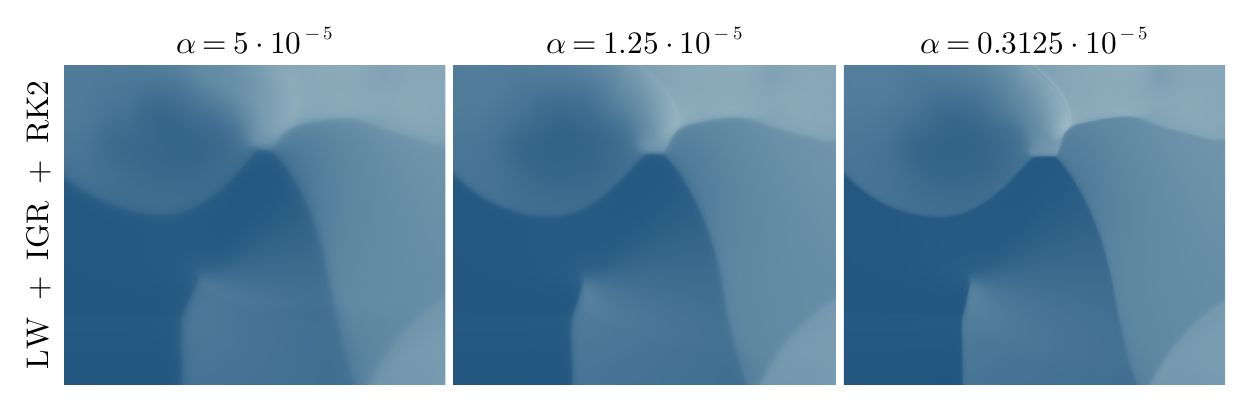}
     \vspace{-8 mm}
     \caption{\textbf{Convergence for $\alpha \rightarrow 0$.} As $\alpha \rightarrow 0$, the density gradients sharpen and the regularized solutions converge to the nominal solution (plots with $\Delta x = \Delta y = 5 \cdot 10^{-4}$).}
     \label{fig:sine_waves_alpha_convergence} 
\end{figure}
LW without IGR broke down with negative densities. Thus, we only show LW + IGR. 
To improve the efficiency of LF and LW + IGR, we integrate them into fourth (for LF) and second-order Runge-Kutta (RK) methods, respectively. 
We discretize the elliptic problem defining $\Sigma$ using a five-point Laplacian stencil and central differences for the $u$-derivatives appearing in the right-hand side. 
At each flux computation, we solve it via Gauss-Seidel iteration \cite{saad2003iterative}, with the previous solution as a warm start. 
Except for \cref{fig:sine_waves_alpha_convergence} that explores impractically large $\alpha$, all results in this section are computed using a single Gauss-Seidel sweep per flux computation.

\subsubsection{Blast wave interaction}
We simulate the interaction of three blast waves in $[0, 0.72] \times [0, 1.2]$, with periodic boundary conditions.
The initial conditions are $u(x, y, 0) \equiv 0$ and $\rho(x, y, 0) = 1 + 0.6 \eta_1(x, y) + 1.2 \eta_2(x, y) + 0.5 \eta_3(x, y)$, where the $\eta_i$ are Gaussian densities with means $((0.2, 0.2), (0.4, 0.7), (0.3, 1.05))$ and standard deviations $(0.05, 0.075, 0.03)$.
\Cref{fig:blast_waves_density} compares LF + RK4 and LW + RK2 on a $432 \times 720$ grid to a reference solution computed with LF + RK4 on a $4320 \times 7200$ grid.
LW + RK2, with $\alpha = 3.6 ((\Delta x)^2 = (\Delta y)^2)$, preserves the fine details at later time steps while the viscous regularization of LF + RK4 smoothes them out.

\subsubsection{Shear shocks}
We investigate the shear-like collision of two shock waves on $[0.0, 1.2] \times [0.0, 1.0]$, with periodic boundary conditions. 
We use initial conditions $u(x, y \leq 0.5, 0) \approx 2.5 \sin(2 \pi (x - 1) / 1.2) + 1$, $u(x, y \geq 0.5, 0) \approx - 1.5 \sin(2 \pi (x - 0.9) / 1.2) - 1$, and $\rho(x, y, 0) \equiv 1$, ensuring a gradual transition at $y=0$ and $y=1$ using a smooth cutoff function. 
The layers $y \leq 0.5$ and $y \geq 0.5$ form shocks moving in opposite directions that collide in the shear layers causing interactions of shock waves with strong vortical flow.
In \cref{fig:sine_waves_density}, we compare LF + RK4 and LW + RK2 on $600 \times 500$ grid to a reference solution computed with LF + RK4 on a $1.2 \cdot 10^4 \times 10^4$ grid.
LW + IGR + RK2, computed with $\alpha = 5 ((\Delta x)^2 = (\Delta y)^2)$, preserves the fine details while the viscous regularization of LF + RK4 smoothes them out.
As shown in \cref{fig:2d_euler_sine_wave_energy}, LW + IGR + RK2 avoids the systematic loss of energy exhibited by LF + RK4.
In \cref{fig:sine_waves_longtime_density} we show that LW + IGR + RK2 preserves sharp density gradients for much longer time horizons. 
In \cref{fig:sine_waves_mesh_convergence}, we provide evidence that refining the mesh for constant $\alpha$ results in convergence to a well-defined solution of \cref{eqn:infgeo_reg_flux_dd}.
Finally, in \cref{fig:sine_waves_alpha_convergence} we illustrate the convergence to the reference solution under decreasing $\alpha$.

\section{Comparison, conclusion, and outlook}
\subsection{Comparison to related work}
\label{sec:comparison}
\begin{figure}
     \begin{tikzpicture}
\begin{groupplot}[
	group style={group size=4 by 1,
	horizontal sep=0.045\textwidth,
    vertical sep=1.05cm,},
	]
	\nextgroupplot[
		title={},
		standard,
 		xlabel={$x$},
		ylabel={$\rho(x, T)$},
		height=0.202\textwidth,
		width=0.202\textwidth,
        ymin=0.0,
        ymax=8.0,
        xmin=0.3,
        xmax=0.7,
		xmajorticks=false,
		ymajorticks=false,
		xtick={},
		ytick={},
		enlarge x limits=0.,
		enlarge y limits=0.,
		legend style={
			at={(2.15, 1.5)},inner sep=4pt,anchor=north,legend columns= 3,cells ={align = right}, draw=none,fill=none}, 
		]	

		\addlegendimage{steelblue, ultra thick}
		\addlegendimage{orange, ultra thick}
		\addlegendimage{rust, ultra thick, dash pattern=on 7pt off 1pt}
		\addlegendimage{empty legend}
		\addlegendimage{seagreen, ultra thick}
		\addlegendimage{darksky, ultra thick , dash pattern=on 6pt off 6pt}

		\legend{Unregularized, LAD with $\frac{\alpha}{(\Delta x)^2} = 2.5 $, IGR with $\frac{\alpha}{(\Delta x)^2} = 2.5$, \phantom{placeholder} , LAD with $\frac{\alpha}{(\Delta x)^2} = 250$, IGR with $\frac{\alpha}{(\Delta x)^2} = 250 $};
		\addplot[orange, ultra thick] table[x index={2},y index={1}] {csv/lad_comparison_oscillation/lad.csv};

	\nextgroupplot[
		standard,
		title={},
		xlabel={$x$},
		ylabel={$\rho(x, T)$},
		height=0.202\textwidth,
		width=0.202\textwidth,
        ymin=0.0,
        ymax=8.0,
        xmin=0.3,
        xmax=0.7,
		xmajorticks=false,
		ymajorticks=false,
		xtick={},
		ytick={},
		enlarge x limits=0.,
		enlarge y limits=0.,
		]	

		\addplot[rust, ultra thick, dash pattern=on 7pt off 1pt, dash phase=5pt] table[x index={2},y index={1}] {csv/lad_comparison_oscillation/igr.csv};

	\nextgroupplot[
		standard,
		title={},
		xlabel={$x$},
		ylabel={$u(x,0)$},
		height=0.202\textwidth,
		width=0.202\textwidth,
        ymin=-0.002,
        ymax=0.002,
        xmin=0.0,
        xmax=1.0,
		xmajorticks=false,
		ymajorticks=false,
		xtick={},
		ytick={},
		enlarge x limits=0.,
		enlarge y limits=0.,
		]	
	
		\addplot[steelblue, very thick, smooth] table[x index={1},y index={0}] {csv/lad_comparison_dissipation/init_u.csv};
		
	\nextgroupplot[
		standard,
		title={},
 		xlabel={$t$},
 		ylabel={Energy},
		height=0.202\textwidth,
		width=0.202\textwidth,
        ymin=0.0,
        ymax=0.00000026,
        xmin=0.0,
        xmax=1.0,
		ymajorticks=false,
		xmajorticks=false,
        yticklabels={{$0$}},
		ytick={0},
		enlarge x limits=0.,
		enlarge y limits=0.,
		]	

		\addplot[steelblue, ultra thick] table[x index={1},y index={0}] {csv/lad_comparison_dissipation/rlw.csv};
		\addplot[rust, ultra thick, , dash pattern=on 7pt off 1pt, dash phase=4.5pt] table[x index={1},y index={0}] {csv/lad_comparison_dissipation/igr_low_alpha.csv};
		\addplot[darksky, ultra thick, dash pattern=on 6pt off 6pt] table[x index={1},y index={0}] {csv/lad_comparison_dissipation/igr_high_alpha.csv};
		\addplot[orange, ultra thick] table[x index={1},y index={0}] {csv/lad_comparison_dissipation/lad_low_alpha.csv};
		\addplot[seagreen, ultra thick] table[x index={1},y index={0}] {csv/lad_comparison_dissipation/lad_high_alpha.csv};

	\end{groupplot}
\end{tikzpicture}
     \vspace{-4mm}
     \caption{\textbf{IGR versus localized artificial diffusivity (LAD).} The nonlocality of IGR reduces oscillations of dispersive schemes like Lax-Wendroff at shocks (first two panels, initial condition is a 1-d Riemann problem). 
     In high-frequency, low-amplitude sound waves (third panel), LAD loses energy (as defined in \cref{fig:2d_euler_sine_wave_energy}) as each downward slope of $u$ incurs dissipation.
     The antidissipation of IGR in upward slopes prevents this, even for 100 times larger $\alpha$. (fourth panel)} 
     \label{fig:IGR_v_LAD}
\end{figure}
The earliest work proposing inviscid regularization in a spirit similar to ours is a series of papers by Bhat and Fetecau \cite{bhat2006hamiltonian,bhat2009riemann,bhat2009regularization}, motivated by ideas of Leray~\cite{leray1934essai}.
In the case of Burgers' equation shows promising results. 
But its generalization to the Euler equation fails to maintain the shock speed of the nominal equation, greatly limiting its usefulness for numerical purposes \cite{bhat2009regularization}.
Recently, \cite{guelmame2022hamiltonian} proposed a family of inviscid and nondispersive regularizations for the unidimensional Euler equation, based on a similar regularization for the shallow water equation \cite{clamond2018non}.
The members of this family can be written in conservation form and thus maintain the correct shock speed. 
They are parametrized by a function $\mathcal{A}: \R \longrightarrow \R$ of the density. 
In the special case of $P \equiv 0$ and $\mathcal{A}(\rho) \coloneqq \rho$, this regularization coincides with \cref{eqn:infgeo_reg_flux_1d} \emph{with $2\alpha$ replaced by $\alpha$} in the right-hand side of the elliptic problem defining $\Sigma$.
This amounts to using the Levi-Civita connection of the Fisher-Rao metric, instead of the dual connection in this work. 
Levi-Civita geodesics move at constant speed, as measured by the Riemannian metric. 
This translates to smooth solutions of the associated PDE conserving the energy given by the squared Riemannian norm. 
Thus, they need to form singularities to model irreversible behavior like the dissipation in shocks.
This is the likely cause for the unphysical cusp singularities in the regularized solutions of \cite{guelmame2022hamiltonian}, as studied by \cite{pu2018weakly,liu2019well}.
Numerous existing works apply viscous regularization adaptively, based on shock indicators \cite{vonneumann1950method,puppo2004numerical,cook2005hyperviscosity,fiorina2007artificial,mani2009suitability,barter2010shock,guermond2011entropy,bruno2022fc,dolejvsi2003some,bhagatwala2009modified,kawai2010assessment}. 
In the unidimensional case, a prototypical form of this ``localized artificial diffusivity'' (LAD) amounts to removing the elliptic operator and using only the negative part of $[\partial_x u]$ in \cref{eqn:infgeo_reg_flux_1d}, yielding
\begin{equation}
     \label{eqn:igr_lad}
     \begin{cases}
          \partial_t
          \begin{pmatrix}
               \rho u \\
               \rho
          \end{pmatrix}
          + \partial_x
          \begin{pmatrix}
               \rho u^2 + P(\rho)  + \Sigma \\ 
               \rho u
          \end{pmatrix}
          = 
          \begin{pmatrix} 
               f\\
               0
          \end{pmatrix} \\
          \rho^{-1} \Sigma \Ccancel[orange]{{- \alpha \partial_x(\rho^{-1} \partial_x \Sigma)}} = 2 \alpha {\color{orange}\big(}[\partial_x u]{\color{orange}\big)_{\operatorname{-}}} [\partial_x u].
     \end{cases}
\end{equation}
Thus, the regularization term $\Sigma$ of \cref{eqn:infgeo_reg_flux_1d} mimics a nonlocal and sign-indefinite diffusion. 
As shown in \cref{fig:IGR_v_LAD}, the nonlocality due to the elliptic problem in IGR reduces spurious oscillations in dispersive numerical schemes for a given parameter $\alpha$.
Where LAD can impose limits on the time step sizes due to the quadratic scaling of parabolic CFL numbers \cite{kawai2010assessment}, our experiments suggest that increasing $\alpha$ in IGR imposes no such restrictions, although it may increase the number of Jacobi or Gauss-Seidel sweeps per time steps needed to compute $\Sigma$.
A striking difference between IGR and LAD is that the IGR ``diffusion tensor'' is sign-indefinite. 
As shown in \cref{fig:IGR_v_LAD}, this overcomes LAD's systematic dissipative bias in high-frequency sound waves that otherwise requires specialized treatment \cite{kawai2010assessment}.
In the multivariate case, IGR avoids energy dissipation due to shear motion, different from diffusive regularization. 
Similar considerations have motivated different variants of adaptive bulk viscosity (ABV) \cite{cook2005hyperviscosity,mani2009suitability} in the numerical solution of Navier Stokes problems. 
Indeed, the variant proposed by \cite{mani2009suitability} (with parameter $r = 0, C = 2\alpha$) amounts to using $\Sigma = 2 \alpha \left(\trace\left(\cnst{D}u\right)\right)_{-} \left(\trace\left(\cnst{D}u\right)\right)$ in IGR. 
But we find that without the $\trace\left([\cnst{D}u]^2\right)$ term, both IGR and ABV are prone to vanishing density in strongly rotational flows. 
The Riemannian geometry induced by the squared $L^2$-distance in the original Euler equation is closely related to the Wasserstein geometry on the space of mass densities \cite{brenier1989least,khesin2021geometric} and the $\operatorname{logdet}$ regularization equals the negative Shannon entropy of the mass density. 
Information geometric regularization compromises between them.
Similarly, \cite{jordan1998variational,amari2018information,amari2019information,peyre2019computational,cui2022time,li2023wasserstein,amari2023information,burger2023covariance} combine optimal transport with entropic terms and information geometry and \cite{keith2023proximal} use entropic regularization to solve explicitly inequality-constrained PDEs.
In providing a geometric approach to modeling the energy dissipation in a shock our work is also related to the geometric treatment of thermodynamic systems in \cite{gay2017lagrangiana,gay2017lagrangianb,gay2018lagrangian,yoshimura2023hamiltonian}.

\subsection{Conclusion}
In this work, we have proposed an information geometric regularization of the barotropic compressible Euler equation and provided a proof of concept for its numerical application. 
To the best of our knowledge, this is the first work that regularizes continuum mechanical equations by modifying the geometry of the manifold of deformation maps. 
It thus opens up numerous directions for future work, including the application to other systems of PDEs, the use of different geometries, as well as the theoretical and numerical study of the resulting equations. 
For the sake of concreteness and conciseness, we have restricted this work to the case of the barotropic Euler equation and the information geometry induced by the logarithmic barrier.
But generalizations to different geometries and equations like Navier Stokes and general Euler equations are straightforward to derive. 
An improved theoretical understanding of our approach will help identify its most promising applications.  
\subsection{Outlook: Toward information geometric mechanics}
\label{sec:igm}
Our work combines the geometry governing the inertial movement of individual particles with the information geometry of distributions of particles.
It thus takes the first step toward developing \emph{information geometric mechanics} that accounts for the geometries not just of fundamental physical laws governing microscopic behavior, but also of the statistical summaries describing their emergent macroscopic behavior.

\section*{Acknowledgments}
The authors gratefully acknowledge support from the Air Force Office of Scientific Research under award number FA9550-23-1-0668 (Information Geometric Regularization for Simulation and Optimization of Supersonic Flow). RC acknowledges support through a PURA travel award from the Georgia Tech Office of Undergraduate Education. 
This work was produced using Microsoft Copilot for autocompletion and spelling/grammar checking.
We thank the anonymous reviewers for their constructive comments that helped us to improve this manuscript.

\bibliographystyle{siamplain}
\bibliography{references}

\end{document}